\documentclass[11pt,leqno]{article}
\usepackage{amssymb,amsmath,amsthm}

\usepackage[all,2cell,ps]{xy}
\usepackage{enumerate}

\newcommand{\symplecto}{\omega}

\UseAllTwocells

\newtheorem{theorem}{Theorem}[section]
\newtheorem{proposition}[theorem]{Proposition}
\newtheorem{lemma}[theorem]{Lemma}
\newtheorem{corollary}[theorem]{Corollary}
\newtheorem{convention}[theorem]{Convention}
\theoremstyle{definition}

\newtheorem{definition}[theorem]{Definition}
\newtheorem{notation}[theorem]{Notation}
\newtheorem{example}[theorem]{Example}

\newtheorem{remark}[theorem]{Remark}

\newtheorem{conjecture}[theorem]{Conjecture}

\newcommand{\C}{{\mathbb{C}}}
\newcommand{\N}{{\mathbb{N}}}
\newcommand{\R}{{\mathbb{R}}}

\newcommand{\Z}{{\mathbb{Z}}}

\newcommand{\kor}{{\mathbb{K}}}
\newcommand{\cor}{{\bf{k}}}
\newcommand{\hcor}{\widehat{\bf{k}}}

\newcommand{\coro}{{\bf{k}}{_0}}
\newcommand{\hcoro}{\widehat{\bf{k}}{_0}}

\def\phi{{\varphi}}
\def\epsilon{\varepsilon}

\def\sha{\mathcal{A}}
\def\A{\mathcal{A}}
\def\shb{\mathcal{B}}
\def\shc{\mathcal{C}}
\def\shd{\mathcal{D}}
\def\she{\mathcal{E}}
\def\shf{\mathcal{F}}

\def\shi{\mathcal{I}}
\def\shj{\mathcal{J}}

\def\shl{\mathcal{L}}
\def\shm{\mathcal{M}}
\def\shn{\mathcal{N}}
\def\sho{\mathcal{O}}

\def\sht{\mathcal{T}}

\def\shw{\mathcal{W}}

\newcommand{\OO}{\sho}

\newcommand{\E}[1][]{\mathcal{E}_{#1}}

\newcommand{\EE}[1][]{\she_{#1}}

\newcommand{\W}[1][]{\mathcal{W}_{#1}}

\newcommand{\HW}[1][]{{\widehat{\mathcal W}}_{#1}}

\newcommand{\Dttau}{\shd_X[\opb{\tau},\tau]}
\newcommand{\Dotau}{\shd_X[\opb{\tau}]}

\newcommand{\WW}[1][]{\shw_{#1}}

\newcommand{\gr}{\mathop{\mathrm{gr}}}

\newcommand{\Lie}[1][]{\operatorname{\mathsf{L}}\def\temp{#1}
\ifx\temp\empty\else^{(#1)}\fi}

\newcommand{\stx}{{\mathfrak{X}}}
\newcommand{\sty}{{\mathfrak{Y}}}

\newcommand{\stkHom}[1][]{\mathfrak{Hom}_{\raise1.5ex\hbox to.1em{}#1}}

\DeclareMathOperator{\coker}{Coker}
\DeclareMathOperator{\im}{Im}

\newcommand{\rmpt}{{\rm pt}}
\newcommand{\rmptt}{{\{\rm pt\}}}

\def\epito{{\twoheadrightarrow}}
\renewcommand{\to}[1][]{\xrightarrow[]{#1}}
\newcommand{\from}[1][]{\xleftarrow[]{#1}}
\newcommand{\isoto}[1][]{\xrightarrow[#1]%
{{\raisebox{-.6ex}[0ex][-.6ex]{$\mspace{1mu}\sim\mspace{2mu}$}}}}

\newcommand{\TopHom}[1][]{\mathrm{TopHom}_{\raise1.5ex\hbox to.1em{}#1}}
\newcommand{\BHom}[1][]{\mathrm{Bhom}_{\raise1.5ex\hbox to.1em{}#1}}

\newcommand{\Hom}[1][]{\mathrm{Hom}_{\raise1.5ex\hbox to.1em{}#1}}
\newcommand{\RHom}[1][]{\mathrm{RHom}_{\raise1.5ex\hbox to.1em{}#1}}
\newcommand{\Ext}[2][]{\mathrm{Ext}_{\raise1.5ex\hbox to.1em{}#1}^{#2}}
\renewcommand{\hom}[1][]{{\mathcal{H}om}_{\raise1.5ex\hbox to.1em{}#1}}
\newcommand{\rhom}[1][]{{R\mathcal{H}om}_{\raise1.5ex\hbox to.1em{}#1}}
\newcommand{\ext}[2][]{{\mathcal{E}xt}_{\raise1.5ex\hbox to.1em{}#1}^{#2}}

\newcommand{\Tens}[1][]{\mathbin{\otimes_{\raise1.5ex\hbox to-.1em{}#1}}}
\newcommand{\LTens}[1][]{\mathbin{\otimes_{\raise1.5ex\hbox to-.1em{}#1}^{L}}}
\newcommand{\Tor}[2][]{\mathrm{Tor}^{\raise1.5ex\hbox to.1em{}#1}_{#2}}
\newcommand{\tens}[1][]{\mathbin{\otimes_{\raise1.5ex\hbox to-.1em{}#1}}}
\newcommand{\ltens}[1][]{\mathbin{\otimes_{\raise1.5ex\hbox to-.1em{}#1}^{L}}}
\newcommand{\tor}[2][]{{\mathcal{T}or}^{\raise1.5ex\hbox to.1em{}#1}_{#2}}
\newcommand{\etens}{\mathbin{\boxtimes}}
\newcommand{\detens}{\underline{\etens}}

\newcommand{\htens}[1][]{{\widehat\tens}_{#1}}

\newcommand{\Endo}[1][]{\mathrm{End}_{\raise1.5ex\hbox to.1em{}#1}}

\newcommand{\Aut}[1][]{\mathrm{Aut}_{\raise1.5ex\hbox to.1em{}#1}}
\newcommand{\sect}{\Gamma}
\newcommand{\rsect}{\mathrm{R}\Gamma}

\newcommand{\RD}{{\rm D}}
\newcommand{\RC}{{\rm C}}
\newcommand{\RK}{{\rm K}}
\newcommand{\Rb}{{\rm b}}
\newcommand{\coh}{{\rm coh}}

\newcommand{\hol}{{\rm hol}}
\newcommand{\rh}{{\rm rh}}
\newcommand{\reg}{{\rm reg}}
\newcommand{\grad}{{\rm grad}}

\newcommand{\wCc}{{{\rm w-}\C\rm -c}}
\newcommand{\Cc}{{\C\rm -c}}

\newcommand{\oim}[1]{{#1}_*}

\newcommand{\reim}[1]{{R#1}_!}
\newcommand{\opb}[1]{#1^{-1}}

\newcommand{\tw}[1]{\widetilde{#1}}

\DeclareMathOperator{\supp}{supp}

\DeclareMathOperator{\chv}{char}

\newcommand{\Fct}{\operatorname{Fct}}
\newcommand{\for}{\mathit{for}}

\def\rop{{\rm op}}
\def\tot{{\rm tot}}

\newcommand{\md[1]}{{\rm Mod}({#1})}
\newcommand{\mdf[1]}{{\rm Mod^f}({#1})}

\newcommand{\mdbc[1]}{{\rm Mod^{bc}}({#1})}
\newcommand{\mdfn[1]}{{\rm Mod^{dfn}}({#1})}

\newcommand{\Mods}{\operatorname{{\mathfrak{Mod}}}}
\newcommand{\Mod}{\mathrm{Mod}}

\newcommand{\indlim}[1][]{\mathop{\varinjlim}\limits_{#1}}

\newcommand{\eqdot}{\mathbin{:=}}

\newcommand{\cl}{\colon}
\newcommand{\scbul}{\,\raise.4ex\hbox{$\scriptscriptstyle\bullet$}\,}
\newcommand{\vvert}{\vert\vert}

\newcommand{\ba}{\begin{array}}
\newcommand{\ea}{\end{array}}

\newcommand{\bnum}{\begin{enumerate}[{\rm(i)}]}
\newcommand{\enum}{\end{enumerate}}
\newcommand{\banum}{\begin{enumerate}[{\rm(a)}]}
\newcommand{\eanum}{\end{enumerate}}

\newcommand{\lp}{{\rm(}}
\newcommand{\rp}{{\rm)}}

\newcommand{\eq}{\begin{eqnarray}}
\newcommand{\eneq}{\end{eqnarray}}
\newcommand{\eqn}{\begin{eqnarray*}}
\newcommand{\eneqn}{\end{eqnarray*}}

\newcommand{\II}[1][]{{\rm Ind}({#1})}
\newcommand{\PP}[1][]{{\rm Pro}({#1})}

\newcommand{\prolim}[1][]{\mathop{\varprojlim}\limits_{#1}}

\newcommand{\proolim}[1][]{\mathop{``{\varprojlim}"}\limits_{#1}}
\newcommand{\inddlim}[1][]{\mathop{``{\varinjlim}"}\limits_{#1}}

\newcommand{\CD}[1][]{\mathcal{C}_{#1}}
\newcommand{\rmw}{{\rm w}}
\newcommand{\set}[2]{\left\{#1 \mathbin{;} #2 \right\}}
\providecommand{\bysame}{\makebox[3em]{\hrulefill}\thinspace}

\numberwithin{equation}{section}
\begin{document}



\title
{Constructibility and duality
for simple holonomic modules on complex symplectic manifolds}
\author{Masaki Kashiwara and Pierre Schapira}
\maketitle
\begin{abstract}
Consider a complex symplectic manifold $\stx$ and  
the algebroid stack $\WW[\stx]$
of deformation-quantization. 
For two regular holonomic  $\WW[\stx]$-modules $\shl_i$ ($i=0,1$) 
supported by smooth Lagrangian  submanifolds,
we prove that the complex $\rhom[{\WW[\stx]}](\shl_1,\shl_0)$ is
perverse
over the  field $\WW[\rmpt]$ 
and dual to the
complex $\rhom[{\WW[\stx]}](\shl_0,\shl_1)$. 
\end{abstract}


\section*{Introduction}
\renewcommand{\thefootnote}{}
\footnote{Mathematics Subject Classification: 46L65, 14A20, 32C38}
Consider a complex symplectic manifold $\stx$. 
A local model for $\stx$ is an
open subset of  the cotangent bundle $T^*X$ to a
complex manifold $X$, and $T^*X$ is 
endowed with the filtered sheaf of rings 
$\shw_{T^*X}$ of deformation-quantization. This sheaf of rings is  
similar to the sheaf $\she_{T^*X}$ of microdifferential operators of
Sato-Kawai-Kashiwara~\cite{S-K-K}, but with an extra
central parameter $\tau$, a substitute to the lack of homogeneity
(see \S~\ref{section:Wmod}). This is an algebra over the field
$\cor=\shw_\rmpt$, a subfield of the field of formal Laurent series 
$\C[[\opb{\tau},\tau]$. 

It would be tempting to glue the locally 
defined sheaves of algebras  $\shw_{T^*X}$ to give rise to a 
sheaf on $\stx$, but the procedure fails, and 
 one is lead to replace  
the notion of a sheaf of algebras by that of an ``algebroid stack''. 
A canonical algebroid stack $\WW[\stx]$
on $\stx$, locally equivalent to $\shw_{T^*X}$,  has been constructed by 
Polesello-Schapira~\cite{P-S} after
Kontsevich~\cite{Ko} had treated the formal case (in the general
setting of Poisson manifolds) by a different method.

In this paper, we study regular holonomic $\WW[\stx]$-modules
supported by 
smooth Lagrangian submanifolds of $\stx$. For example,  a regular holonomic
module along the zero-section $T^*_XX$ of $T^*X$ is 
locally isomorphic to a finite sum of copies of 
the sheaf $\sho_X^\tau$ whose sections are
series 
$\sum_{-\infty<j \leq m}f_j\tau^j$ ($m\in\Z$), where the $f_j$'s
are sections of $\sho_X$ and  satisfy certain
growth conditions.

Denote by $\RD^\Rb_{\rh}(\WW[\stx])$ the full subcategory  
of the bounded derived
category of $\WW[\stx]$-modules on $\stx$ consisting of objects with
regular holonomic cohomologies.
The main theorem of this paper asserts that if $\shl_i$ ($i=0,1$) 
are objects of $\RD^\Rb_{\rh}(\WW[\stx])$ 
supported by smooth Lagrangian manifolds 
$\Lambda_i$ and if one sets 
$F\eqdot\rhom[{\WW[\stx]}](\shl_1,\shl_0)$, then 
$F$ is $\C$-constructible over 
the  field $\cor$, its microsupport is contained in the
normal cone $C(\Lambda_0,\Lambda_1)$ and $F$ is dual over $\cor$ to the object 
$\rhom[{\WW[\stx]}](\shl_0,\shl_1)$. If  $\shl_0$ and   $\shl_1$ are
concentrated in degree $0$, then $F$ is perverse.
We make the conjecture that the
hypothesis on the smoothness of the $\Lambda_i$'s  may be removed. 

The strategy of our proof is as follows.

First, assuming only that $\shl_0$ and $\shl_1$ are coherent, 
 we construct the canonical morphism
\eq\label{eq:dual00}
&&\rhom[{\WW[\stx]}](\shl_1,\shl_0)\to
\RD_\stx\bigl(\rhom[{\WW[\stx]}](\shl_0,\shl_1\,[\dim_\C\stx])\bigr),
\eneq
where $\RD_\stx$ is the duality functor for sheaves of
$\cor_\stx$-modules.  By using a kind of Serre's duality 
for the sheaf $\sho_X^\tau$, we prove that \eqref{eq:dual00} is an isomorphism 
as soon as $\rhom[{\WW[\stx]}](\shl_0,\shl_1)$ is constructible. For
that purpose, 
we need to develop some functional analysis over nuclear algebras, in the line 
of Houzel \cite{Ho}.

In order to prove the constructibility result for
$F=\rhom[{\WW[\stx]}](\shl_1,\shl_0)$, we  may assume that 
$\stx=T^*X$,  $\Lambda_0=T^*_XX$, $\shl_0=\sho_X^\tau$ and $\Lambda_1$
is the graph of the
differential of a holomorphic function $\phi$ defined on $X$. 
Consider the sheaf of rings $\shd_X[\opb{\tau}]= \shd_X\tens\C[\opb{\tau}]$.  
We construct a coherent $\shd_X[\opb{\tau}]$-module $\shm$
which generates $\shl_1$ and, setting
$F_0\eqdot\rhom[{\shd_X[\opb{\tau}]}](\shm,\sho_X^\tau(0))$, we prove that 
the microsupport of $F_0$   is 
a closed complex analytic Lagrangian subset of $T^*X$ contained in
$C(\Lambda_0,\Lambda_1)$. Using  a deformation argument as in  \cite{K1}  
and some functional analysis (extracted from \cite{Ho}) over 
the base ring $\coro=\W[\rmpt](0)$,  
we deduce that the fibers of the cohomology of $F_0$ are finitely
generated, and the result follows from the isomorphism
$F\simeq F_0\tens[\coro]\cor$.

In this paper, we also consider a complex compact contact manifold $\sty$ and 
the algebroid
stack $\EE[\sty]$ of microdifferential operators on it. 
Denote by $\RD^\Rb_{\rh}(\EE[\sty])$ the $\C$-triangulated category
consisting of $\EE[\sty]$-modules with regular holonomic cohomology. 
Using results of Kashiwara-Kawai~\cite{K-K}, 
we prove that this category has finite $\Ext[]{}$,
admits a Serre functor in the sense of Bondal-Kapranov~\cite{B-K} and
this functor is nothing but a shift by $\dim_\C\,\sty +1$. In
other words,  $\RD^\Rb_{\rh}(\EE[\sty])$ is a Calabi-Yau category. 
Similar results hold over $\cor$ for a complex compact symplectic manifold $\stx$.

\section{$\shw$-modules on $T^*X$}\label{section:Wmod}

Let $X$ be a complex manifold, $\pi\cl T^*X\to X$ its cotangent
bundle. The homogeneous symplectic manifold $T^*X$ 
is endowed with the   $\C^\times$--conic $\Z$-filtered sheaf of 
rings $\E[T^*X]$  of finite--order  microdifferential operators,
and its subring $\E[T^*X](0)$ of operators of order $\leq 0$ 
constructed in \cite{S-K-K} (with other notations). 
We assume that the reader is familiar with the 
theory of $\she$-modules, referring to \cite{K3} or \cite{S} for an
exposition. 

On the symplectic manifold $T^*X$ there exists  another (no more conic)
sheaf of rings $\shw_{T^*X}$,  called ring of  deformation-quantization  
by many authors (see Remark~\ref{re:formal} below). 
The study of the relations between the ring $\E[T^*X]$  on a 
complex homogeneous symplectic 
 manifold and the sheaf $\shw_{T^*X}$ 
on a complex symplectic manifold is systematically 
performed in~\cite{P-S}, where it is shown in particular 
how quantized symplectic transformations act 
on $\shw_{T^*X}$-modules. We follow here their presentation. 

Let $t\in\C$ be the coordinate and set
$$
\E[T^*(X\times\C),\hat t]=\set{P\in\E[T^*(X\times\C)]}{
[P,\partial_t] = 0}.
$$
Set $T^*_{\tau\neq 0}(X\times\C) = \set{(x,t;\xi,\tau)\in T^*(X\times\C)}%
{\tau\neq 0}$, and
consider the map
\eq\label{eq:projrho}
\rho\cl T^*_{\tau\neq 0}(X\times\C) \to T^*X
\eneq
given in local coordinates by $\rho(x,t;\xi,\tau) = (x;\xi/\tau)$.
The ring $\W[T^*X]$   on $T^*X$ is given by
\eqn
&&\W[T^*X]\eqdot \oim\rho(\E[X\times\C, \hat t]
\vert_{T_{\tau\neq 0}^*(X\times\C)}).
\eneqn
In a local symplectic coordinate system $(x;u)$ on $T^*X$, 
a section $P$ of $\W[T^*X]$ on an open subset $U$ 
is written as a formal series, called its
total symbol:
\eq\label{eq:totsymb}
&&\sigma_{\tot}(P) = \sum_{-\infty<j\leq m}p_j(x;u)\tau^j, \,m\in\Z\, 
\quad p_j\in\OO_{T^*X}(U),
\eneq
with the condition  
\eq\label{eq:defW}
&&\left\{  \parbox{300 pt}{
for any compact subset $K$ of $U$ there exists a positive 
constant $C_K$ such that 
$\sup\limits_{K}\vert p_{j}\vert \leq C_K^{-j}(-j)!$ for all $j<0$.
}\right. 
\eneq
The product is given by the Leibniz rule.
If $Q$ is an operator of total symbol $\sigma_{\tot}(Q)$, then 
$$\sigma_{\tot}(P\circ Q)
=\sum_{\alpha\in\N^n} \dfrac{\tau^{-\vert\alpha\vert}}{\alpha !} 
\partial^{\alpha}_u\sigma_{\tot}(P)\partial^{\alpha}_x\sigma_{\tot}(Q).
$$ 

Denote by $\W[T^*X](m)$ the subsheaf of $\W[T^*X]$ consisting of sections 
$P$ whose total symbol $\sigma_{\tot}(P)$ satisfies: $p_j $ is zero for $j>m$. 
The ring $\W[T^*X]$ is $\Z$-filtered by the $\W[T^*X](m)$'s, and $\W[T^*X](0)$
is a subring. If $P\in\W[T^*X](m)$ and $P\notin\W[T^*X](m-1)$, $P$ is said
of order $m$. Hence $0$ has
order $-\infty$.  
Then 
there is a   well-defined principal symbol morphism, which does not depend on
the local coordinate system on $X$:
\eq
&&\sigma_m\cl \W[T^*X](m)\to\sho_{T^*X}\cdot\tau^m.
\eneq
 If $P\in\W[T^*X]$ has order $m$ on a connected open subset of $T^*X$,
$\sigma_m(P)$ is called the principal symbol of $P$ and is denoted
$\sigma(P)$. Note that  a section $P$ in $\W[T^*X]$ is 
invertible on an open subset $U$ 
of $T^*X$ if and only if its principal symbol is nowhere vanishing on $U$.

The principal symbol map
induces an isomorphism of graded algebras
$$\gr\W[T^*X] \isoto \OO_{T^*X}[\tau,\tau^{-1}].$$
We set 
\eqn
&&\cor \eqdot \W[\rmptt],\quad \cor(j)\eqdot \W[\rmptt](j)\, (j\in\Z),   
\quad\coro \eqdot \cor(0).
\eneqn
Hence, $\cor$ is a  field and an element  $a\in\cor$ is 
written as a formal series
\eq\label{eq:cor}
a = \sum_{j\leq m}a_j\tau^j, \quad a_j\in\C,\quad  m\in\Z
\eneq
satisfying \eqref{eq:defW}, that is,
\eq\label{eq:defWb}
&&\left\{  \parbox{300 pt}{
 there exists a positive constant $C$ such that 
$\vert a_{j}\vert \leq C^{-j}(-j)!$ for all $j<0$.
}\right. 
\eneq
Note that $\cor_0$ is a discrete valuation ring.

\begin{convention}
In the sequel, we shall identify $\partial_t$ and its total symbol
$\tau$. Hence, we consider $\tau$ as a section of $\W[T^*X]$. 
\end{convention}
Note that 
\begin{itemize}
\item
$\W[T^*X]$ is  flat over $\W[T^*X](0)$
and in particular $\cor$ is  flat over $\coro$,
\item
$\coro$ is faithfully flat over $\C[\opb\tau]$, 
\item
 the sheaves of rings $\W[T^*X]$ and $\W[T^*X](0)$ are right
and left Noetherian (see \cite[Appendix]{K3}) 
and in particular coherent, 
\item
if $\shm$ is a coherent $\shw_{T^*X}$-module, its support is a closed
complex analytic involutive (by Gabber's theorem)  subset of $T^*X$,
\item
$\W[T^*X]$ is a $\cor$-algebra and 
$\W[T^*X](0)$ is a $\coro$-algebra,
\item
there are natural monomorphisms of sheaves of $\C$-algebras 
\eq\label{eq:DinW}
&&\opb{\pi}\shd_X \hookrightarrow \E[T^*X]\hookrightarrow\W[T^*X].
\eneq
\end{itemize}
On an affine chart,  
morphism  \eqref{eq:DinW} is described on symbols as follows. 
To a section of $\E[T^*X]$ of total symbol
$\sum_{-\infty<j\leq m}p_j(x;\xi)$ (the $p_j$'s are homogeneous in $\xi$ of degree $j$)
one associates the section of $\W[T^*X]$ of total symbol
$\sum_{-\infty<j\leq m}p_j(x;u)\tau^j$, with $u=\opb{\tau}\xi$. 

\begin{remark}\label{re:formal} 
(i) Many authors consider the filtered ring of formal operators,
defined by
\eqn
\HW[T^*X](m)=\prolim[j\leq m]\W[T^*X](m)/\W[T^*X](j)\quad m\in\Z,
&&\HW[T^*X]=\bigcup_m\HW[T^*X](m).
\eneqn
Then
\eqn
&&\hcor \eqdot \HW[\rmptt]\simeq\C[[\opb\tau,\tau], 
\quad\hcoro \eqdot \HW[\rmptt](0)\simeq\C[[\opb\tau]].
\eneqn
Note that 
$\HW[T^*X]$ is faithfully flat over $\W[T^*X]$ and $\HW[T^*X]$ is  flat over
$\HW[T^*X](0)$.

\noindent
(ii) Many authors  also prefer to use the symbol $\hbar=\opb{\tau}$ instead
of $\tau$.
\end{remark}

\begin{notation}\label{not:shotau}
Let $\shi_X$
be the left ideal of $\W[T^*X]$   generated by the vector fields on $X$.
We set
\eqn
&&\OO^\tau_X =\W[T^*X]/\shi_X,
    \quad\OO^\tau_X(m)=\W[T^*X](m)/(\shi_X\cap \W[T^*X](m))\, (m\in\Z).
\eneqn
\end{notation}

Note that $\OO^\tau_X$ is a coherent 
$\W[T^*X]$-module supported by the zero-section $T^*_XX$.  A section $f(x,\tau)$
of this module may be written as a series:
\eq
&& f(x,\tau)=\sum_{-\infty<j\leq m}f_j(x)\tau^j,\quad m\in\Z,
\eneq
the $f_j$'s satisfying the condition \eqref{eq:defW}.

Also note that $\sho_X$ is a direct summand  of $\sho_X^\tau$ as a sheaf. 

\begin{lemma}\label{le:shotau=sho}
After identifying $X$ and $X\times\{0\}\subset X\times\C$, 
there is an isomorphism of sheaves of $\C$-vector spaces 
\lp not of algebras\rp\,
\eq\label{eq:shotau=sho}
&& \sho_X^\tau(0)\simeq \sho_{X\times\C}\vert_{X\times\{0\}}.
\eneq
\end{lemma}
\begin{proof}
By its construction, $\sho_X^\tau(0)$ is isomorphic to the sheaf of 
holomorphic microfunctions $\shc_{X\times\{0\}\vert X\times\C }(0)$ of
\cite{S-K-K},
and this last sheaf is isomorphic to
$\sho_{X\times\C}\vert_{X\times\{0\}}$ by loc.\ cit.
\end{proof}
Note that isomorphism \eqref{eq:shotau=sho} corresponds to the map
\eqn
&&\sho_X^\tau(0)\ni\sum_{j\leq 0}f_j\tau^{j}
\mapsto \sum_{j\geq 0}f_{-j}\dfrac{t^{j}}{j!}\in
\sho_{X\times\C}\vert_{X\times\{0\}}.
\eneqn

We shall use the following sheaves of rings.
We set 
\eq\label{defDtau}
\Dotau \eqdot \shd_X\tens\C[\opb{\tau}],
&&\Dttau \eqdot \shd_X\tens\C[\tau,\opb{\tau}].
\eneq
Note that 
$\OO_X^\tau(0)$  is a left $\Dotau$-module and 
$\OO_X^\tau$  a left $\Dttau$-module. 

\begin{lemma}\label{le:dtau1}
\bnum
\item
$\Dotau$ is a right and left Noetherian sheaf  of\/ $\C$-algebras,
\item
if $\shn\subset\shm$ are two coherent $\Dotau$-modules and $\shm_0$ a
coherent $\shd_X$-submodule of $\shm$, then 
$\shn\cap\shm_0$ is $\shd_X$-coherent,
\item
if $\shm$ is  a coherent $\Dotau$-module and $\shm_0$ a
$\shd_X$-submodule of $\shm$ of finite type, then $\shm_0$ is
$\shd_X$-coherent.
\enum
\end{lemma}
\begin{proof}
(i) follows from \cite[Th. A.~3]{K3}.

\noindent
(ii) Set $\sha\eqdot\Dotau$ and let $\sha_n\subset\sha$ be the
$\shd_X$-submodule consisting of sections of order $\leq n$ with respect
to $\opb{\tau}$. 
We may reduce to the case where 
$\shm=\sha^N$. By \cite[Th.~A.29, Lem.~A21]{K3}, $\shn\cap\sha_n$ is
$\shd_X$-coherent. Since $\{\shn\cap\sha_n\cap\shm_0\}_n$ is an
increasing sequence of coherent $\shd_X$-submodules of $\shm_0$ and 
this last module is finitely generated, the sequence is stationary.

\noindent
(iii) Let us keep the notation in (ii). Here again, we 
may reduce to the case where $\shm=\sha^N$. Then $\shm_0\cap\sha_n$
is $\shd_X$-coherent, and the proof goes as in (ii). 
\end{proof}

\begin{lemma}\label{le:dtau2}
\bnum
\item
$\Dttau$ is flat over $\Dotau$,
\item
if $\shi$ is a finitely generated left ideal of $\Dttau$, then
$\shi\cap\Dotau$ is a locally finitely generated left ideal of $\Dotau$.
\enum
\end{lemma}
\begin{proof}
(i) $\C[\tau,\opb{\tau}]$ is flat over $\C[\opb{\tau}]$.

\noindent
(ii) Let $\shi_0\subset\shi$ be a coherent $\Dotau$-module which
generates $\shi$. Set $\shj_n\eqdot\tau^n\shi_0\cap\Dotau$. Then 
$\shi\cap\Dotau=\bigcup_n\shj_n$ and this increasing sequence 
of coherent ideals of $\Dotau$ is locally stationary. 
\end{proof}

\begin{lemma}\label{le:dtau3}
$\W[T^*X]$ is flat over $\opb{\pi}\Dttau$. 
\end{lemma}
\begin{proof}
One proves that $\W[T^*X]$ is flat over $\opb{\pi}\Dotau$ exactly as one
proves that the ring of microdifferential operators $\E[T^*X]$
is flat over $\opb{\pi}\shd_X$. 
Since we have $\Dttau\tens[\Dotau]\shm\simeq\shm$ 
for any $\Dttau$-module $\shm$, 
the result follows.
\end{proof}

\section{Regular holonomic $\W$-modules}

The following definition  adapts a classical
definition of \cite{K-O} to $\W[T^*X]$-modules  (see also \cite{D-S}).

\begin{definition}\label{def:simpleWmod}
Let $\Lambda$ be a smooth Lagrangian submanifold of $T^*X$. 
\banum
\item
Let $\shl(0)$ be a coherent $\W[T^*X](0)$-module supported by $\Lambda$. 
One says that $\shl(0)$ is regular (resp.~simple) along
$\Lambda$ if
$\shl(0)/\shl(-1)$ is a coherent $\OO_\Lambda$-module 
(resp.\ an invertible $\OO_\Lambda$-module). 
Here,  $\shl(-1)=\W[T^*X](-1)\shl(0)$.
\item
Let $\shl$ be a coherent $\W[T^*X]$-module supported by $\Lambda$.  
One says that $\shl$ is regular (resp.~simple) along
$\Lambda$ if there exists locally a coherent $\W[T^*X](0)$-submodule $\shl(0)$
of $\shl$ such that $\shl(0)$  generates $\shl$ over $\W[T^*X]$ and is 
regular (resp.~simple) along $\Lambda$.
\eanum
\end{definition}
Note that in (b), $\shl(0)/\shl(-1)$ is necessarily a locally free 
$\OO_\Lambda$-module.
\begin{example}
The sheaf $\OO^\tau_X$ is a simple $\W[T^*X]$-module along the
zero-section $T^*_XX$.
\end{example}

The following result is easily proved (see \cite{D-S}).
\begin{proposition}\label{pr:locSimp}
Let $\Lambda$ be a smooth Lagrangian submanifold of $T^*X$  and let
$\shm$ be a coherent $\W[T^*X]$-module. 
\bnum
\item
If $\shm$ is regular, then it is locally a finite direct sum of simple modules.
\item
Any two $\W[T^*X]$-modules simple along 
$\Lambda$ are locally
isomorphic. In particular, any simple module along $T^*_XX$ is 
locally isomorphic to $\OO^\tau_X$.
\item
If $\shm$, $\shn$ are simple $\W[T^*X]$-modules along
$\Lambda$, 
then
$\rhom[{\W[T^*X]}](\shm, \shn)$ is concentrated in degree $0$ and is 
a $\cor$-local system of rank one on 
$\Lambda$.
\enum
\end{proposition}

\begin{definition}\label{def:regholWmod}
Let $\shm$ be a coherent $\W[T^*X]$-module and let $\Lambda$ denote its 
support (a closed $\C$-analytic subset of $T^*X$). 
\bnum
\item
We say that $\shm$ is holonomic if $\Lambda$ is Lagrangian.
\item
Assume $\shm$ is holonomic. We say that $\shm$ is regular
holonomic if there is an open subset $U\subset\stx$ such that 
$U\cap\Lambda$ is a dense subset of the regular locus $\Lambda_\reg$
of $\Lambda$ and $\shm\vert_U$ is regular
along $U\cap\Lambda$. 
\enum
\end{definition}
In other words, a holonomic module $\shm$ is regular if it is so at
the generic points of its support.
Note that when $\Lambda$ is smooth, Definition~\ref{def:regholWmod}
is compatible with  Definition~\ref{def:simpleWmod}.
This follows from Gabber's theorem. Indeed, 
we have the following theorem, analogous of \cite[Th.~7.34]{K3}:
\begin{theorem}
Let $U$ be an open subset of $T^*X$, $\shm$ a coherent
$\shw_{T^*X}\vert_U$-module and $\shn\subset\shm$ a
sub-$\shw_{T^*X}(0)$-module. Assume that $\shn$ is a 
small filtrant inductive limit 
of coherent $\shw_{T^*X}(0)$-modules. Let $V$ be the set of $p\in U$ in a
neighborhood of which $\shn$ is $\shw_{T^*X}(0)$-coherent. Then 
$U\setminus  V$ is an analytic involutive subset of $U$. 
\end{theorem}
The fact that regularity of holonomic $\shw$-modules is a generic
property follows as in loc.\ cit.~Prop.~8.28. 

\section{$\W[]$-modules on a complex symplectic manifold}

We refer to \cite{Ko} for the definition of an algebroid stack and
to \cite{D-P} for a more systematic study. 

Let $\kor$ be a
commutative unital ring and let $X$ be a topological space.
We denote by 
$\md[\kor_X]$ the abelian category of sheaves of 
$\kor$-modules and by $\RD^\Rb(\kor_X)$ its bounded derived category.

If $A$ is a $\kor$-algebra, we denote by $A^+$ the category with 
one object and having $A$ as endomorphisms of this object. 
If $\sha$ is a sheaf of $\kor$-algebras on 
$X$, we denote by $\sha^+$ the $\kor$-linear stack associated with
the prestack $U\mapsto \sha(U)^+$ ($U$ open in $X$) and call it
the $\kor$-algebroid stack associated with $\sha$.  
It is equivalent to the stack 
of right $\A$-modules locally isomorphic to $\A$,
and $\A$-linear homomorphisms.

The projective cotangent bundle $P^*Y$
to a complex manifold $Y$ is endowed with the sheaf of rings 
$\she_{P^*Y}$ of microdifferential operators. (This sheaf is
 the direct image of the sheaf
$\she_{T^*Y}$ of \S~\ref{section:Wmod} by the map 
$T^*Y\setminus T^*_YY\to P^*Y$.) 
A complex contact manifold $\sty$ is locally isomorphic
to an open subset of a projective cotangent bundle $P^*Y$ and on 
such a contact manifold, a canonical 
algebroid stack $\EE[\sty]$ locally equivalent to the 
stack associated with the sheaf of rings 
$\she_{P^*Y}$  has been constructed in~\cite{K2}.

This construction has been adapted to the symplectic case by  
\cite{P-S}.
A complex symplectic manifold $\stx$ is locally isomorphic
to the cotangent bundle $T^*X$ to a complex manifold $X$ and  
 a canonical algebroid stack $\WW[\stx]$ locally equivalent to the stack 
associated with  the sheaf of
 rings $\W[T^*X]$ of \S~\ref{section:Wmod} is constructed in loc.\ cit., 
after Kontsevich~\cite{Ko} had  treated 
 the general case of complex  Poisson manifolds 
in the formal  setting by a different approach. 

Denote by $\stx^a$ the complex manifold $\stx$ endowed with 
the symplectic form  $-\omega$, where  $\omega$ is the symplectic form 
on $\stx$. There is a natural equivalence of
algebroid stacks $\WW[\stx^a]\simeq(\WW[\stx])^\rop$. 

Let $\stx$ and $\sty$ be two complex symplectic manifolds. 
There exist a natural 
$\cor$-algebroid stack $\WW[\stx]\etens\WW[\sty]$
on $\stx\times\sty$ and a natural functor of $\cor$-algebroid stacks
$\WW[\stx]\etens\WW[\sty]\to \WW[\stx\times\sty]$ which locally
corresponds to the morphism of sheaves of rings 
$\W[T^*X]\etens\W[T^*Y]\to \W[T^*(X\times Y)]$.

One sets
\eqn
&&\md[{\WW[\stx]}]=\Fct_\cor(\WW[\stx],\Mods(\cor_\stx)),
\eneqn
where 
$\Fct_\cor(\scbul,\scbul)$ denotes the $\cor$-linear category
of $\cor$-linear functors of stacks and $\Mods(\cor_\stx)$ is the
 stack of sheaves of $\cor$-modules on $\stx$. We denote by
 $\Mods(\WW[\stx])$ the stack on $\stx$ given by  
$U\mapsto \md[{\WW[\stx]\vert_U}]$.

Then $\md[{\WW[\stx]}]$ is a Grothendieck abelian category.
We denote by $\RD^\Rb(\WW[\stx])$ its bounded derived category and
call an object of this derived category a
$\WW[\stx]$-module on $\stx$, for short. 
Objects of  $\md[{\WW[\stx]}]$ are described with some details 
in \cite{D-S}. 

We denote by $\hom[{\WW[\stx]}]$ the 
 hom-functor  of the stack  $\Mods(\WW[\stx])$, a functor 
from 
$\md[{\WW[\stx]}]^\rop\times\md[{\WW[\stx]}]$ to $\md[\cor_\stx]$.
The object $\rhom[{\WW[\stx]}](\shm,\shn)$ is thus well defined in
$\RD^\Rb(\cor_{\stx})$ for two objects $\shm$ and $\shn$ of $\RD^\Rb(\WW[\stx])$. 

The natural functor $\WW[\stx]\etens\WW[\sty]\to \WW[\stx\times\sty]$
defines a functor of stacks
\eq\label{eq:tensetens}
&& \for\cl \Mods(\WW[\stx\times\sty])\to \Mods(\WW[\stx]\etens\WW[\sty])
\eneq
and this last functor admits an adjoint (since it does locally). 
For $\shm\in\RD^\Rb(\WW[\stx])$ and  
$\shn\in\RD^\Rb(\WW[\sty])$, it thus exists a canonically defined object
$\shm\detens\shn\in \RD^\Rb(\WW[\stx\times\sty])$
such that, locally,  
$\shm\detens\shn\simeq\WW[T^*(X\times Y)]\etens_{(\WW[T^*X]\etens\WW[T^*Y])}(\shm\etens\shn)$
in $\RD^\Rb(\WW[T^*(X\times Y)])$.

Being local, the notions of coherent or holonomic, regular holonomic
or simple 
object of $\md[{\WW[\stx]}]$ make sense. We denote by
\begin{itemize}
\item
$\RD^\Rb_{\coh}(\WW[\stx])$ the full triangulated subcategory 
 of $\RD^\Rb(\WW[\stx])$ consisting
of objects with coherent cohomologies,
\item 
$\RD^\Rb_{\hol}(\WW[\stx])$ the full triangulated subcategory 
 of $\RD^\Rb_{\coh}(\WW[\stx])$ consisting 
of objects with holonomic cohomologies, or equivalently, of objects with
Lagrangian supports in $\stx$,
\item 
$\RD^\Rb_{\rh}(\WW[\stx])$ the full triangulated 
subcategory of $\RD^\Rb_{\hol}(\WW[\stx])$ consisting of objects with
regular holonomic cohomologies.
\end{itemize}
The support of an object $\shm$ of $\RD^\Rb_{\coh}(\WW[\stx])$ is
 denoted by $\supp(\shm)$ and is also called its characteristic
 variety.
This is a  closed complex analytic involutive subset of $\stx$. 

In the sequel, we shall denote by $\Delta_\stx$ the diagonal of
$\stx\times\stx^a$ and identify it with $\stx$ by the first projection.

The next lemma follows from general considerations on stacks and
its verification is left to the reader. 
\begin{lemma}
There exists a canonical simple $\W[\stx\times\stx^a]$-module $\CD[\Delta_\stx]$
supported by the diagonal $\Delta_\stx$ such that 
if $U$ is open in $\stx$ and isomorphic to an open
subset $V$ of a cotangent bundle $T^*X$, then $\CD[\Delta_\stx]\vert_U$
is isomorphic to $\W[T^*X]\vert_V$ as a $\W[T^*X]\tens(\W[T^*X])^\rop$-module.
\end{lemma}
\begin{definition}
Let $\shm\in\RD^\Rb(\W[\stx])$. 
Its dual  $\RD'_{\rmw}\shm\in\RD^\Rb(\W[\stx^a])$ 
is given by
\eq\label{def:dual1}
&&\RD'_{\rmw}\shm\eqdot \rhom[{\W[\stx]}](\shm,\CD[\Delta_\stx]).
\eneq
\end{definition}
Let $\Lambda$ be a smooth Lagrangian submanifold of $\stx$.
Consider the cohomology class ${\bf c}_{1/2} \in H^2(X;\C_\Lambda^\times)$
defined as  ${\bf c}_{1/2}=\beta(\frac12\alpha([\Omega_\Lambda]))$
where $[\Omega_\Lambda]\in H^1(X;\sho_\Lambda^\times)$ 
is the class of the line bundle $\Omega_\Lambda$
and $\alpha,\beta$ are  the morphisms of the exact sequence
\eqn
&&H^1(X;\sho_\Lambda^\times) \to[\alpha] H^1(X;d\sho_\Lambda)
\to[\beta]H^2(X;\C_\Lambda^\times).
\eneqn
The main theorem of \cite{D-S} asserts that 
 simple $\W[\stx]$-modules along $\Lambda$ are in one-to-one
correspondence with  twisted local systems of rank one on 
$\Lambda$ with twist ${\bf c}_{1/2}$.

\begin{proposition}\label{pro:etensandhom}
Let $\shm,\shn$ be two objects of 
$\RD^\Rb_{\coh}(\WW[\stx])$. There is a natural isomorphism
in $\RD^\Rb(\cor_\stx)$:
\eqn
&&\rhom[{\WW[\stx]}](\shm,\shn)\simeq\rhom[{\WW[\stx\times\stx^a]}]
(\shm\detens\RD'_{\rmw}\shn,\CD[\Delta_\stx]).
\eneqn
\end{proposition}
\begin{proof}
We have the isomorphism in $\RD^\Rb(\WW[\stx])$:
\eqn
&&\shn\simeq \rhom[{\WW[\stx^a]}](\RD'_{\rmw}\shn,\CD[\Delta_\stx]),
\eneqn
from which we deduce
\eqn
\rhom[{\WW[\stx\times\stx^a]}]
(\shm\detens\RD'_{\rmw}\shn,\CD[\Delta_\stx]) &\simeq& 
\rhom[{\WW[\stx]}]
(\shm,\rhom[{\WW[\stx^a]}](\RD'_{\rmw}\shn,\CD[\Delta_\stx]))\\
&\simeq& \rhom[{\WW[\stx]}](\shm,\shn).
\eneqn
\end{proof}

\section{Functional analysis I}

In this section, we will use techniques elaborated by 
Houzel \cite{Ho} and will follow his terminology. 
(See also Kiel-Verdier \cite{K-V} for related results.)

We call a bornological convex $\C$-vector space (resp.\ $\C$-algebra),
a $bc$-space (resp.\ $bc$-algebra) and we
denote by $\mdbc[\C]$ the category of $bc$-spaces and bounded linear
maps. This additive category admits small inductive and projective limits, but
is not abelian.
 
Let $A$ be a $bc$-algebra. 
We denote by $\mdbc[A]$ the additive category 
of bornological $A$-modules and bounded $A$-linear maps.
For $E,F\in\mdbc[A]$, we set:
\eqn
&& \BHom[A](E,F)=\Hom[{\mdbc[A]}](E,F),\\
&& E^\vee=\BHom[A](E, A).
\eneqn

Let $E\in\mdbc[A]$ and let $B\subset E$ be a convex circled bounded subset of
$E$. For $x\in E$, one sets 
\eqn
&& \vvert x\vvert_B=\mathop{\inf}\limits_{x\in cB, c\in \C}\vert c\vert.
\eneqn

For $u\in \BHom[A](E,F)$,  $B$ bounded in $E$ and  $B'$ convex circled 
bounded in $F$,  one sets 
\eqn
&& \vvert u\vvert_{BB'}=\mathop{\sup}\limits_{x\in B}\vvert u(x)\vvert_{B'}.
\eneqn

One says that a sequence $\{u_n\}_n$ in $\BHom[A](E,F)$ is bounded if
for any bounded subset $B\subset E$ there exists a convex 
circled bounded subset
$B'\subset F$ such that $\mathop{\sup}\limits_{n}\vvert u_n\vvert_{BB'}<\infty$.

One says that $u\in \BHom[A](E,F)$ is $A$-nuclear if there exist
a bounded sequence $\{y_n\}_n$ in $F$, a bounded sequence $\{u_n\}_n$
in $\BHom[A](E,A)$ and a summable sequence $\{\lambda_n\}_n$ in
$\R_{\geq 0}$ such that
\eqn
&& u(x)=\sum_n\lambda_n\,u_n(x)\,y_n\text{ for all }x\in E.
\eneqn

For $E,F\in\mdbc[A]$, there is
a natural structure of $bc$-space on $E\tens[A]F$ and one
denotes by $E\htens[A]F$ the completion of $E\tens[A]F$.
Assuming $F$ is complete, there is a natural linear map 
\eqn
&& E^\vee\htens[A]F \to \BHom[A](E,F).
\eneqn
An  element $u\in \BHom[A](E,F)$ 
is $A$-nuclear if and only if it is in the image of $E^\vee\htens[A]F $. 

Note that, if  $u\cl E\to F$ is a nuclear $\C$-linear map, 
$u\htens \, 1\cl E\htens\,A\to  F\htens\,A$ is $A$-nuclear.

Recall that a  $\C$-vector space $E$ is called  
a DFN-space if it is an inductive limit $\{(E_n,u_n)\}_{n\in\N}$ of Banach
spaces such that the maps $u_n\cl E_n\to E_{n+1}$ are
$\C$-nuclear and injective. Note that any bounded subset of $E$ is
contained in $E_n$ for some $n$.

In the sequel, we will consider a DFN-algebra $A$ and the
full subcategory $\mdfn[A]$ of $\mdbc[A]$ consisting of
DFN-spaces. Note that any epimorphism $u\cl E\to F$ in $\mdbc[A]$ is
semi-strict, that is, any bounded sequence in $F$ is the image by $u$
of a bounded sequence in $E$. Also note that for 
$E$ and $F$ in $\mdfn[A]$, $\BHom[A](E,F)$ is the subspace of 
$\Hom[A](E,F)$ consisting of continuous maps.
Since the category $\mdfn[A]$ is not abelian, we introduce the
following definition, referring to \cite{Scn0} for a more systematic
treatment of homological algebra in terms of quasi-abelian categories.

\begin{definition}
\bnum
\item
Let $A$ be a DFN-algebra.
A complex $0\to E'\to E\to E''\to 0$ in  $\mdfn[A]$ is called a short exact sequence 
if it is an exact sequence in $\md[A]$.
\item
Let $A$ and $B$ be two DFN-algebra. An additive functor from 
$\mdfn[A]$ to $\mdfn[B]$ is called exact if it sends short exact
sequences to short exact sequences.
\enum
\end{definition}

The following result is well-known and follows from \cite{Gr}.
\begin{proposition}\label{pro:htensexactindfn}
Let $A$ be a DFN-algebra. The functor
$\scbul\htens\,A\cl\mdfn[\C]\to\mdfn[A]$ is exact.
\end{proposition}

Recall that a $bc$-algebra $A$ is multiplicatively convex if for any
bounded set $B\subset A$, there exist a constant $c>0$ and 
a convex circled bounded
set $B'$ such that $B\subset c\cdot B'$ and $B'\cdot B'\subset B'$.

As usual, for an additive category $\shc$, 
we denote by $\RC^\Rb(\shc)$ the category of bounded
complexes in $\shc$ and, 
for $a\leq b$ in $\Z$, by $\RC^{[a,b]}(\shc)$ the full
subcategory consisting of complexes concentrated in degrees $j\in [a,b]$. 
We denote by $\RK^\Rb(\shc)$ the homotopy category associated
with $\RC^\Rb(\shc)$. Finally, we denote by
$\II[\shc]$ and $\PP[\shc]$ the categories of ind-objects and
pro-objects of $\shc$, respectively. 

\begin{theorem}\label{th:functan1}
Let $A$ be a multiplicatively convex  DFN-algebra and assume that $A$
is a Noetherian ring \lp when forgetting the topology\rp. 
Consider an
inductive system $\{(E^\bullet_n,u^\bullet_n)\}_{n\in\N}$ in 
$\RC^{[a,b]}(\mdfn[A])$ for $a\leq b\in\Z$.
Assume:
\bnum
\item
$u^\bullet_n\cl E_n^\bullet\to E_{n+1}^\bullet$ is a quasi-isomorphism
for all $n\geq 0$,
\item
$u_n^j\cl E_n^j\to E_{n+1}^j$ is $A$-nuclear for all $j\in\Z$ and all
$n\geq 0$.
\enum
Then $H^j(E_n^\bullet)$ is finitely generated over $A$ for all $j\in\Z$ and
all $n\geq 0$.  
\end{theorem}
This  is a particular case of \cite[\S~6~Th.~1, Prop.~A.1]{Ho}.

\begin{theorem}\label{th:functan2}
Let $A$ be a  DFN-algebra, and 
consider an inductive system  $\{E^\bullet_n,u^\bullet_n\}_{n\in\N}$ 
in $\RC^{[a,b]}(\mdfn[A])$ for $a\leq b$ in $\Z$.
Assume
\bnum
\item
for each $i\in \Z$ and $n\in \N$, the map $u^i_n\cl E^i_n\to E^i_{n+1}$ is
$A$-nuclear,
\item
for each $i\in \Z$, $\inddlim[n] H^i(E^\bullet_n)\simeq 0$
in $\II[{\Mod(\C)}]$. 
\enum
Then $\inddlim[n] E^\bullet_n\simeq 0$ in $\II[{\RK^{[a,b]}(\mdfn[A])}]$.
\end{theorem}
First, we need a lemma.

\begin{lemma}\label{le:functan3}
Consider the solid diagram in $\mdfn[A]$:
\eqn
\xymatrix{
 E\ar[r]^-u\ar@{.>}[d]_-v&F\ar[d]^-{v'}\\
 E'\ar[r]_-{u'}&F'.
}\eneqn
Assume that $u$ is $A$-nuclear and $\im v'\subset \im u'$. Then there
exists a morphism $v\cl E\to E'$ making the whole diagram commutative.
\end{lemma}
\begin{proof}
The morphism $w\cl E'\times_{F'}F\to F$ is well defined in the
category $\mdfn[A]$ and is surjective by the hypothesis. We get a
diagram 
\eqn
\xymatrix{
 &E\ar[d]^-u\\
 E'\times_{F'}F'\ar@{->>}[r]_-{w}&F, 
}\eneqn
In this situation, the nuclear map $u$ factors
through $E'\times_{F'}F'$ by \cite[\S4~Cor.~2]{Ho}.
\end{proof}

\begin{proof}[Proof of Theorem~\ref{th:functan2}]
We may assume that $H^i(E^\bullet_n)\to H^i(E^\bullet_{n+1})$ is
the zero morphism for all $i\in\Z$ and all $n\in\N$. 
Consider the solid diagram
\eqn
\xymatrix{
 E^b_{n-1}\ar[r]^-{u^b_{n-1}}\ar@{.>}[d]_-{k^b_n}&E^b_n\ar[d]^-{u^b_n}\\
 E^{b-1}_{n+1}\ar[r]_-{d^{b-1}_{n+1}}&F^b_{n+1}.
}\eneqn
Since $H^b(E^\bullet_p)\simeq \coker d^{b-1}_p$ for all $p$, 
$\im u^b_n\subset \im d^{b-1}_{n+1}$. Moreover, $u^b_{n-1}$ is
$A$-nuclear by the hypothesis. Therefore we may apply Lemma~\ref{le:functan3}
and we obtain a map $k^b_n\cl E^b_{n-1}\to E^{b-1}_{n+1}$ making the
whole diagram commutative. Set $v^i_n=u^i_n\circ u^i_{n-1}$ and 
$h^b_n=d^{b-1}_{n+1}\circ k^b_n$. Consider the diagram
\eqn
\xymatrix{
{\cdots}\ar[rr]&&E^{b-2}_{n-1}\ar[rr]^-{d^{b-2}_{n-1}}\ar[d]^-{v^{b-2}_n}
                &&E^{b-1}_{n-1}\ar[rr]^-{d^{b-1}_{n-1}}\ar[d]^-{v^{b-1}_n}
                  &&E^{b}_{n-1}\ar[rr]\ar[d]^-{v^{b}_n}\ar[lld]_-{h^{b}_n}
                    &&0\\
{\cdots}\ar[rr]&&E^{b-2}_{n+1}\ar[rr]^-{d^{b-2}_{n+1}}
                &&E^{b-1}_{n+1}\ar[rr]^-{d^{b-1}_{n+1}}
                  &&E^{b}_{n+1}\ar[rr]
                    &&0.
}\eneqn
The morphisms $v^i_n$'s define a morphism of complex 
$v_n\cl E^\bullet_{n-1}\to E^\bullet_{n+1}$. We define
$h_n^i\cl  E^i_{n-1}\to E^i_{n+1}$ by setting $h^i_n=0$ for $i\neq b$. 
Now denote by $\sigma^{\leq b-1}E^\bullet_n$ the stupid
truncated complex obtained by replacing $E^b_n$ with $0$. The morphism
\eqn
&&v_n- h_n\circ d_{n-1}-d_{n+1}\circ h_n\cl E^\bullet_{n-1}\to E^\bullet_{n+1}
\eneqn
factorizes through $\sigma^{\leq b-1}E^\bullet_n$. 
Hence, we get an isomorphism
\eqn
&&\inddlim[n]E^\bullet_n\isoto\inddlim[n]\sigma^{\leq b-1}E^\bullet_{n}
\eneqn
 in $\II[{\RK^{[a,b]}(\mdfn[A])}]$.
By repeating this
argument, we find the isomorphism 
$\inddlim[n]E^\bullet_n\simeq \inddlim[n]\sigma^{\leq
  b-p}E^\bullet_{n}$ for any $p\in\N$. 
This completes the proof.
\end{proof}

\begin{theorem}\label{th:functan3}
Let $A$ be a  DFN-field, let $a\leq b$ in $\Z$, 
consider an inductive system  $\{E^\bullet_n,u^\bullet_n\}_{n\in\N}$ in
$\RC^{[a,b]}(\mdfn[A])$
and set $F^\bullet_n=(E^\bullet_n)^\vee=\BHom[A](E^\bullet_n,A)$.
Assume
\bnum
\item
for each $i\in \Z$ and $n\in \N$, the map $u^i_n\cl E^i_n\to E^i_{n+1}$ is
$A$-nuclear,
\item
for each $i\in \Z$, $\inddlim[n] H^i(E^\bullet_n)$ belongs to
$\mdf[A]$ \lp the category of finite-dimensional $A$-vector spaces\rp. 
\enum
Then we have the isomorphism
\eq\label{eq:isofunctan3}
&&\proolim[n] F^\bullet_n\simeq \proolim[n]\Hom[A](E^\bullet_n,A)
\text{ in }\PP[{\RK^\Rb(\md[A])}].
\eneq 
In particular, for each $i\in\Z$, $\proolim[n] H^{-i}(F^\bullet_n)$
belongs to $\mdf[A]$ and is
dual of $\inddlim[n] H^i(E^\bullet_n)$.
\end{theorem}
\begin{proof}
Recall (see \cite[Exe~15.1]{K-S2} that for an abelian category $\shc$ and 
an inductive system $\{X_j\}_{j\in J}$ in $\RD^{[a,b]}(\shc)$ indexed by a small
filtrant category $J$, if the object $\inddlim[j]H^i(X_j)$ of
$\II[\shc]$ is representable for all 
$i\in\Z$, then $\inddlim[j]X_j\in \II[\RD^\Rb(\shc)]$ is representable.

Applying this result to our situation, we find that the object 
$\inddlim[n]E^\bullet_n$ of $\II[{\RD^\Rb(\md[A])}]$ is 
representable in $\RD^\Rb(\md[A])$. 

Denote by $L^\bullet$ the complex given by  $L^i=\inddlim[n] H^i(E^\bullet_n)$ 
and zero differentials. Since $A$ is  a field, there
exists an isomorphism
$L^\bullet\isoto \inddlim[n]E^\bullet_n$ in $\RD^\Rb(\md[A])$,
hence  a quasi-isomorphism
$u\cl L^\bullet\to \inddlim[n]E^\bullet_n$
in $\RC^\Rb(\md[A])$. There  exists
$n\in\N$ such that $u$ factorizes through 
$L^\bullet\to E^\bullet_n$ for some $n$ and we may
assume $n=0$. Since $L^\bullet$ belongs to $\RC^\Rb(\mdf[A])$, 
$L^\bullet\to E^\bullet_0$ is well defined in $\RC^\Rb(\mdfn[A])$.
For any $n$, let $u_n\cl L^\bullet\to E^\bullet_n$ be the induced morphism.
Let $G^\bullet_n$ be the mapping cone of $u_n$. 
The morphism $u_n$ is embedded in a distinguished triangle in $\RK^\Rb(\mdfn[A])$
\eqn
&& L^\bullet\to[u_n] E^\bullet_n\to G^\bullet_n\to[+1].
\eneqn  
By Theorem~\ref{th:functan2}, $\inddlim[n]H^i(G^\bullet_n)\simeq 0$  
implies $\inddlim[n]G^\bullet_n\simeq 0$ in $\II[{\RK^\Rb(\mdfn[A])}]$. 
Hence,  for any $K\in\RK^\Rb(\mdfn[A])$, the morphism
\eqn
\Hom[{\RK^\Rb(\mdfn[A])}](K,L^\bullet)\to \indlim[n]
\Hom[{\RK^\Rb(\mdfn[A])}](K,E^\bullet_n)
\eneqn
is an isomorphism. By the Yoneda lemma, we have thus obtained the
isomorphism 
\eqn
&&L^\bullet\isoto\inddlim[n] E^\bullet_n \text{ in }\II[{\RK^\Rb(\mdfn[A])}].
\eneqn
Since $L^\bullet$ belongs to  $\RK^\Rb(\mdf[A])$,
$(L^\bullet)^\vee\simeq \Hom[A](L^\bullet,A)$
and we obtain 
\eqn
&&\proolim[n]F^\bullet_n\simeq (L^\bullet)^\vee \simeq
\Hom[A](L^\bullet,A)
\simeq \proolim[n]\Hom[A](E_n^\bullet,A)
\eneqn
in $\PP[{\RK^\Rb(\md[A])}]$.
The isomorphisms for the cohomologies follow
 since $\proolim$ commutes with $H^{-i}(\scbul)$ and
 $\Hom[A](\scbul,A)$. 
\end{proof}

A similar result to Theorem~\ref{th:functan3} 
holds for projective system. 

\begin{theorem}\label{th:functan3bis}
Let $A$ be a  DFN-field, let $a\leq b$ in $\Z$, 
consider a projective system  $\{F^\bullet_n,v^\bullet_n\}_{n\in\N}$ in
$\RC^{[a,b]}(\mdfn[A])$
and set $E^\bullet_n=\BHom[A](F^\bullet_n,A)$.
Assume
\bnum
\item
for each $i\in \Z$ and $n\in \N$, the map $v^i_{n}\cl F^i_{n+1}\to F^i_{n}$ is
$A$-nuclear,
\item
for each $i\in \Z$, $\proolim[n] H^i(F^\bullet_n)$ belongs to
$\mdf[A]$ \lp the category of finite-dimensional $A$-vector spaces\rp. 
\enum
Then we have the isomorphism
\eq\label{eq:isofunctan3bis}
&&\inddlim[n] E^\bullet_n\simeq \inddlim[n]\Hom[A](F^\bullet_n,A)
\text{ in }\II[{\RK^\Rb(\md[A])}].
\eneq 
In particular, for each $i\in\Z$, $\inddlim[n] H^{-i}(E^\bullet_n)$
belongs to $\mdf[A]$ and is
dual of $\proolim[n] H^i(F^\bullet_n)$.
\end{theorem}
The proof being similar to the one of Theorem~\ref{th:functan3}, we
shall not repeat it. 

In the course of \S~\ref{section:funct2} below, we shall also need the 
next lemma.

\begin{lemma}\label{le:factorizdualmor}
Let $A$ be a DFN-algebra and let $u\cl E_0\to E_1$ be a $\C$-nuclear
map of DFN-spaces. 
Recall that $(\scbul)^\vee=\BHom[A](\scbul,A)$ and set  
$(\scbul)^\star=\BHom[\C](\scbul,\C)$.
Then the solid commutative diagram below may be
completed with the dotted arrow as a commutative diagram:
\eqn
\xymatrix{
E_1^\star\htens A\ar[r]^-{u^\star\htens A}\ar[d]&E_0^\star\htens A\ar[d]\\
(E_1\htens A)^\vee\ar[r]_-{(u\htens A)^\vee}\ar@{.>}[ru]&(E_0\htens A)^\vee.
}\eneqn
\end{lemma}
\begin{proof}
Consider the commutative diagram
\eqn
\xymatrix{
E_0^\star\htens E_1\ar[d]\ar[r]&{\BHom[\C](E_0,E_1)}\ar[d]\\
 {\BHom[A]((E_1\htens A)^\vee,E_0^\star\htens A)}\ar[r]
          &{\BHom[A]((E_1\htens A)^\vee,(E_0\htens A)^\vee)}.
}\eneqn
Since $u$ is nuclear, it is the image of an element of
$E_0^\star\htens E_1$.
\end{proof}

\section{Functional analysis II}\label{section:funct2}

This section will provide the framework for apply 
Theorems~\ref{th:functan1} and~\ref{th:functan3}.

Here, $X$ will denote a complex manifold.
For a locally free $\sho_X$-module $\shf$ of finite
rank, we set
\eqn
&&\shf^\tau\eqdot \shf\tens[\sho_X]\sho_X^\tau
\quad \shf^\tau(0)\eqdot \shf\tens[\sho_X]\sho_X^\tau(0).
\eneqn

\begin{lemma}\label{le:functan0}
\bnum
\item
The $\C$-algebra $\coro$ is  a multiplicatively convex  DFN-algebra,
\item
the $\C$-algebra $\cor$ is  a   DFN-algebra,
\item
the functor $\scbul\htens\,\coro\cl\mdfn[\C]\to\mdfn[\coro]$ is
well-defined and exact,
\item
the functor $\scbul\htens\,\cor\cl\mdfn[\C]\to\mdfn[\cor]$ is
well-defined and exact, and is isomorphic to the functor 
$(\scbul\htens\,\coro)\tens[\coro]\,\cor$.
\enum
\end{lemma}
\begin{proof}
(i) Define the subalgebra $\coro(r)$ of $\coro$ by 
\eq\label{eq:defWd}
&& \left\{\parbox{300 pt}{
$u=\sum_{j\leq 0}a_j\tau^j$ belongs to $\coro(r)$ if and only if \\
$\vert u\vert_r \eqdot\sum_{j\leq 0} \dfrac{r^{-j}}{(-j)!}\,\vert a_{j}\vert <\infty$.
} \right.
\eneq
Then, for $u$,$v$ in $\coro(r)$, we have 
\eqn
&& \vert u\cdot v\vert_r\leq \vert u\vert_r \cdot\vert v\vert_r.
\eneqn 
Hence, $(\coro(r),\vert\scbul\vert_r)$ is a Banach algebra and 
$\coro$ is multiplicatively convex since it is 
the inductive limit of 
the $\coro(r)$'s. Moreover,  $\coro$ is a DFN-space because the linear maps 
 $\coro(r)\to \coro(r')$ are nuclear for $0<r'<r$.

\noindent
(ii)--(iv) are clear.
\end{proof}
Note that $\cor$ is not multiplicatively convex.

Let $M$ be a real analytic manifold, $X$ a complexification of
$M$. We denote as usual by $\sha_M$ the sheaf on $M$ of real analytic functions,
that is, $\sha_M=\sho_X\vert_M$. Recall that, for $K$ compact in 
$M$, $\sect(K;\sha_M)$ is a DFN-space.
We set
\eq
&&\sha_M^\tau=\sho_X^\tau\vert_M.
\eneq
\begin{lemma}\label{le:functanreal1}
Let $M$ be a real analytic manifold and $K$ a compact subset of $M$. Then  
\bnum
\item 
the sheaf $\sha_M^\tau$ is $\sect(K;\scbul)$-acyclic,
\item
$\sect(K;\sha_M^\tau)\simeq \sect(K;\sha_M)\htens\,\cor $,
\item
the same result holds with $\sha_M^\tau$ and $\cor$ replaced with 
$\sha_M^\tau(0)$ and $\coro$, respectively.
\enum
\end{lemma}
\begin{proof}
Applying Lemma \ref{le:shotau=sho}, 
 we have isomorphisms for each holomorphically convex compact subset $K$ of $X$:
\eqn
\sect(K;\sho^\tau_X(0))&\simeq&\sect(K\times\{0\};\sho_{X\times\C})\\
&\simeq&\sect(K;\sho_X)\htens\, \sho_{\C,0}\\
&\simeq&\sect(K;\sho_X)\htens\,\coro.
\eneqn
This proves (i)--(ii) for $\sha_M^\tau(0)$ and $\coro$. The other case 
follows since 
$\sha_M^\tau\simeq \sha_M^\tau(0)\tens[\coro]\cor$.
\end{proof}

Let us denote, as usual, by $\shd b_M$ the sheaf 
of Schwartz's distribution on $M$. Recall that 
$\sect_c(M;\shd b_M)$ is a DFN-space.

\begin{lemma}\label{le:functanreal2}
Let $M$ be a real analytic manifold. 
There is a unique \lp up to unique isomorphism\rp\, sheaf of
$\cor$-modules $\shd b_M^\tau$ on $M$ which is soft and satisfies
\eqn
&&\sect_c(U;\shd b_M^\tau)\simeq \sect_c(U;\shd b_M)\htens\cor
\eneqn
for each open subset $U$ of $M$.
The same result holds with $\cor$ replaced with $\coro$. In this case, 
we denote by $\shd b_M^\tau(0)$ the sheaf of $\coro$-modules so
obtained. 
\end{lemma}
\begin{proof}
For two open subsets $U_0$ and $U_1$, the sequence 
\eqn
&&0\to \sect_c(U_0\cap U_1;\shd b_M)\htens\,\cor\to
(\sect_c(U_0;\shd b_M)\htens\,\cor)\oplus (\sect_c(U_1;\shd b_M)\htens\,\cor)\\
&&\hspace{6.7cm}\to\sect_c(U_0\cup U_1;\shd b_M)\htens\,\cor\to 0
\eneqn
is exact, and similarly with $\cor$ replaced with $\coro$. 
The results then easily follow. 
\end{proof}

Denote by $\overline X$ the complex conjugate manifold to $X$ and
by $X_\R$ the real underlying manifold, identified with the diagonal of 
$X\times \overline X$. We shall write for short 
$\sha_X^\tau$ and $\shd b_X^\tau$ instead of $\sha_{X_\R}^\tau$ and $\shd b_{X_\R}^\tau$, 
respectively.
We set $\sha_X^{(p,q),\tau}=\sha_X^{(p,q)}\tens[\sha_X]\sha_X^\tau$,
and similarly with  $\shd b_X$ instead of $\sha_X$.

Consider the Dolbeault-Grothendieck complexes of  sheaves of $\cor$-modules
\eq\label{eq:DGa}
&&0\to \sha_X^{(0,0),\tau}\to[\overline\partial]\cdots
\to[\overline\partial] \sha_X^{(0,d),\tau}\to 0,\\
&&0\to \shd b_X^{(0,0),\tau}\to[\overline\partial]\cdots
\to[\overline\partial] \shd b_X^{(0,d),\tau}\to 0.\label{eq:DGdb}
\eneq

\begin{lemma}\label{le:functanreal3}
Both complexes \eqref{eq:DGa} and \eqref{eq:DGdb} are qis to 
$\sho_X^\tau$. The same result holds when replacing 
$\sha_X^{\tau}$, $\shd b_X^{\tau}$ and $\sho_X^\tau$ with 
$\sha_X^{\tau}(0)$, $\shd b_X^{\tau}(0)$ and $\sho_X^\tau(0)$,
respectively.  
\end{lemma} 
The easy proof is left to the reader.

\begin{lemma}\label{le:functan2}
Let $X$ be a complex Stein manifold, $K$ a holomorphically convex compact subset and 
$\shf$  a locally free $\sho_X$-module of finite rank. Then
\bnum
\item 
one has an isomorphism 
$\sect(K;\shf^\tau)\simeq \sect(K;\shf)\htens\,\cor$, 
\item
the $\C$-vector space $\sect(K;\shf^\tau)$ is naturally endowed with a
topology of DFN-space,
\item
$\rsect(K;\shf^\tau)$ is concentrated in degree $0$,
\item
for $K_0$ and $K_1$ two compact subsets of $X$ such that 
$K_0$ is contained in the interior of  $K_1$, the morphism 
$\sect(K_1;\shf^\tau)\to\sect(K_0;\shf^\tau)$  is
$\cor$-nuclear,
\item
the same results hold with $\shf^\tau$ and $\cor$ replaced with 
$\shf^\tau(0)$ and $\coro$, respectively.
\enum
\end{lemma}
\begin{proof}
(i)--(iii) By the hypothesis, the sequence
\eq\label{eq:DGaK}
&&0\to \sect(K;\sho_X)\to 
 \sect(K;\sha_X^{(0,0)})\to[\overline\partial]\cdots
\to[\overline\partial] \sect(K;\sha_X^{(0,d)})\to 0
\eneq
is exact. It remains exact after applying the functor
$\scbul\htens\,\cor$. The result will follow when comparing the
sequence so obtained with the complex  
\eq\label{eq:DGaKtau}
&&0\to \sect(K;\sho_X^\tau)\to 
 \sect(K;\sha_X^{(0,0),\tau})\to[\overline\partial]\cdots
\to[\overline\partial] \sect(K;\sha_X^{(0,d),\tau})\to 0,
\eneq
The case of $\shf$ is treated similarly, replacing 
$\sha_X$ with $\shf\tens[\sho_X]\sha_X$.  

\noindent
(iv) follows from (i) and the corresponding result for $\sho_X$. 

\noindent
(iv) The proof for 
$\shf^\tau(0)$ and $\coro$ is similar.
\end{proof}

\begin{lemma}\label{le:vanishshotau1}
Let $X$ be a complex manifold of complex dimension $d$ and let $\shf$
be a locally free $\sho_X$-module of finite rank. 
Assume that $X$ is Stein. Then 
\bnum
\item
one has the isomorphism 
$H^d_c(X;\shf^\tau)\simeq H^d_c(X;\shf)\htens \cor$,
\item
the $\C$-vector space $ H^d_c(X;\shf^\tau)$ is naturally endowed with a
topology of DFN-space,
\item
 $\rsect_c(X;\shf_X^\tau)$ is concentrated in degree $d$, 
\item
for $U_0$ and $U_1$ two Stein open subset of $X$ with
$U_0\subset\subset U_1$,
the map $H^d_c(U_0;\shf_X^\tau)\to H^d_c(U_1;\shf_X^\tau)$ is $\cor$-nuclear,
\item
the same results hold with $\shf^\tau$ and $\cor$ replaced with 
$\shf^\tau(0)$ and $\coro$, respectively.
\enum
\end{lemma}
\begin{proof}
The proof is similar to that of Lemma~\ref{le:functan2}.
By the hypothesis, the sequence
\eq\label{eq:DGdbc}
&&0\to
 \sect_c(X;\shd b_X^{(0,0)})\to[\overline\partial]\cdots
\to[\overline\partial] \sect_c(X;\shd b_X^{(0,d)})\to
H^d_c(X;\sho_X)\to 0
\eneq
is exact and $H^d_c(X;\sho_X)$ is a DFN-space. 
This sequence will remain exact after applying the functor
$\scbul\htens\cor$. The result will follow when comparing the
obtained sequence with 
\eq\label{eq:DGdbctau}
&&0\to \sect_c(X;\shd b_X^{(0,0),\tau})\to[\overline\partial]\cdots
\to[\overline\partial] \sect_c(X;\shd b_X^{(0,d),\tau})\to
H^d_c(X;\sho_X^\tau)\to 0.
\eneq
The case of $\shf$ is treated similarly, replacing 
$\shd b_X$ with $\shf\tens[\sho_X]\shd b_X$.  
\end{proof}

\begin{proposition}\label{pro:directshot}
Let $f\cl X\to Y$ be a morphism of complex manifolds of complex
dimension $d_X$ and $d_Y$, respectively.
There is a  $\cor$-linear morphism 
\eq\label{eq:tauintgr}
&&\int_f\cl\reim{f}\Omega_X^\tau\,[d_X]\to \Omega_Y^\tau\,[d_Y]
\eneq
functorial with respect to $f$, which sends
$\reim{f}\Omega_X^\tau(0)\,[d_X]$ to $\Omega_Y^\tau(0)\,[d_Y]$
and which 
induces the classical integration morphism 
$\int_f\cl\reim{f}\Omega_X\,[d_X]\to \Omega_Y\,[d_Y]$ when 
identifying $\Omega_X$ and $\Omega_Y$
with a direct summand of $\Omega^\tau_X$  and $\Omega_Y^\tau$, respectively.

In particular, there is a  $\cor$-linear morphism 
\eq\label{eq:tauresidu1}
&&\int_X\cl H^{d_X}_c(X;\Omega_X^\tau)\to \cor
\eneq
which 
induces the classical residues morphism 
$H^{d_X}_c(X;\Omega_X)\to \C$ when identifying 
$H^{d_X}_c(X;\Omega_X)$ \lp resp.\ $\C$\rp\,
with a direct summand of $H^{d_X}_c(X;\Omega^\tau_X)$ \lp resp.\ $\cor$\rp.
\end{proposition}
In the particular case where  $X$ is Stein and $Y=\rmpt$, this follows 
easily from Lemma~\ref{le:vanishshotau1}. Since we shall not use 
the general case, we leave the proof to the reader.

\begin{proposition}\label{pro:serreOtau}
Let $X$ be a Stein complex manifold of complex dimension $d$ and let 
$K$ be a holomorphically convex compact
subset of $X$. Then the  pairing for a Stein open subset $U$ of $X$
\eq\label{eq:inesohtau}
H^d_c(U;\Omega_X^\tau)\times \sect(U;\sho_X^\tau)\to\cor,
&&(fdx,g)\mapsto \int_Ugfdx, 
\eneq
defines the isomorphisms
$\proolim[U\supset K]\BHom[\cor](\sect(\overline U;\sho_X^\tau),\cor)
\simeq\proolim[U\supset K]H^d_c(U;\Omega_X^\tau)$ 
and 
$\inddlim[U\supset K]\BHom[\cor](H^d_c(U;\Omega_X^\tau),\cor)
\simeq\inddlim[U\supset K]\sect(\overline U;\sho_X^\tau)$.
Here, $U$ ranges over the family of Stein open neighborhoods of $K$. 
\end{proposition}
\begin{proof}
We shall prove the first isomorphism, the other case being similar.

By Lemmas~\ref{le:functan2}, \ref{le:vanishshotau1} 
and~\ref{le:factorizdualmor}, we have:
\eqn
\proolim[U\supset K]\BHom[\cor](\sect(\overline U;\sho_X^\tau),\cor)
&\simeq& \proolim[U\supset K](\sho_X(\overline U)\htens \,\cor)^\vee\\
&\simeq& \proolim[U\supset K](\sho_X(\overline U))^\star\htens \,\cor\\
&\simeq& \proolim[U\supset K](\sho_X(U))^\star\htens \,\cor\\
&\simeq& \proolim[U\supset K](H^d_c(U;\Omega_X))\htens \,\cor\\
&\simeq& \proolim[U\supset K]H^d_c(U;\Omega_X^\tau).
\eneqn
\end{proof}

In the proof of Corollary~\ref{cor:main1}, we shall need the
following result.
\begin{lemma}\label{le:vanisectYotau}
Let $Y$ be a smooth closed submanifold of codimension $l$ of $X$. Then 
$H^j(\rsect_Y(\sho_X^\tau))$ and $H^j(\rsect_Y(\sho_X^\tau(0)))$ vanish for $j<l$.
\end{lemma}
\begin{proof}
Since the problem is local, we may assume that $X=Y\times \C^l$ 
where $Y\subset X$ is identified with $Y\times\{0\}$. 

Let $U_n$ be  the open ball of $\C^l$ centered at $0$ with radius
$1/n$.  Using the Mittag-Leffler theorem (see \cite[Prop.~2.7.1]{K-S1}), it
 is enough to prove that for any holomorphically convex
compact subset $K$ of $Y$
\eqn
&&H^j_c(K\times U_n;\sho_X^\tau)=0\text{ for }j<l,
\eneqn
and similarly with $\sho_X^\tau(0)$.
It is enough to prove the result for the sheaf
$\sho_X^\tau(0)$. Then we may replace  $\sho_X^\tau(0)$ with 
the sheaf $\sho_{X\times\C}\vert_{X\times\{0\}}$ and we are
reduced to the well-known result
\eqn
&&H^j_c(K\times U_n\times\{0\};\sho_{X\times\C})=0\text{ for }j<l.
\eneqn
\end{proof}


\begin{remark}
Related results have been obtained, in a slightly different framework, 
in \cite{Pr-S}
\end{remark}

\section{Duality for $\shw$-modules}

We mainly follow the notations of \cite{K-S1}.
Let $X$ be a real manifold and $\kor$ a field.
For $F\in\RD^\Rb(\kor_X)$, we denote by $SS(F)$ its microsupport, a
closed $\R^+$-conic subset of $T^*X$. Recall that this set 
is {\em involutive} (see loc.\ cit.\ Def.~6.5.1).

We denote by $\RD_X$ the duality functor:
\eqn
&&\RD_X\cl (\RD^\Rb(\kor_X))^\rop\to \RD^\Rb(\kor_X),\quad 
F\mapsto \rhom[\kor_X](F,\omega_X),
\eneqn
where $\omega_X$ is the dualizing complex.

Assume now that $X$ is a complex manifold. We denote by $\dim_\C X$ its complex
dimension. We identify the orientation sheaf on $X$ with the
constant sheaf, and the dualizing complex 
$\omega_X$ with $\kor_X\,[2\dim_\C X]$. 

Recall that
an object $F\in\RD^\Rb(\kor_X)$ is
weakly-$\C$-constructible if there exists a complex analytic stratification
$X=\bigsqcup_{\alpha\in A}X_\alpha$ such that
$H^j(F)\vert_{X_\alpha}$ is  locally constant for all $j\in\Z$ and all
$\alpha\in A$. The object $F$ is $\C$-constructible
if moreover $H^j(F)_x$ is finite-dimensional for all
$x\in X$ and all $j\in\Z$.
We denote by  $\RD^\Rb_\wCc(\kor_X)$ 
the full subcategory of  $\RD^\Rb(\kor_X)$ 
consisting of weakly-$\C$-constructible  objects 
and by $\RD^\Rb_\Cc(\kor_X)$ the full subcategory consisting of
$\C$-constructible objects.

Recall (\cite{K-S1}) that $F\in\RD^\Rb(\kor_X)$ is
weakly-$\C$-constructible if and only if its microsupport is a closed
$\C$-conic complex analytic Lagrangian subset of $T^*X$ or, 
equivalently, if it is contained in a closed
$\C$-conic complex analytic isotropic subset of $T^*X$.

{}From now on, our base field is $\cor$.

\begin{theorem}\label{th:duality1}
Let $\stx$ be a complex symplectic manifold and let 
 $\shl_0$ and $\shl_1$ be two objects of $\RD^\Rb_{\coh}(\WW[\stx])$. 
\bnum
\item
There is a natural morphism
\eq\label{eq:dualmorph}
\rhom[{\WW[\stx]}](\shl_1,\shl_0)\to
\RD_\stx\bigl(\rhom[{\WW[\stx]}](\shl_0,\shl_1\,[\dim_\C\stx])\bigr).
\eneq
\item
Assume that $\rhom[{\WW[\stx]}](\shl_0,\shl_1)$
belongs to $\RD^\Rb_{\Cc}(\cor_\stx)$.
Then the morphism~\eqref{eq:dualmorph} is an isomorphism.
In particular, $\rhom[{\WW[\stx]}](\shl_1,\shl_0)$
belongs to $\RD^\Rb_{\Cc}(\cor_\stx)$.
\enum
\end{theorem}

\begin{proof}
By Proposition~\ref{pro:etensandhom}, we may assume that $\shl_0$ 
is a  simple module along a smooth Lagrangian
manifold $\Lambda_0$. In this case, 
$\cor_{\Lambda_0}\to \rhom[{\WW[\stx]}](\shl_0,\shl_0)$ 
is an isomorphism.

\vspace{0.2cm}
\noindent
(i) The natural $\cor$-linear morphism
\eqn
&&\rhom[{\WW[\stx]}](\shl_1,\shl_0)\tens[\cor] \rhom[{\WW[\stx]}](\shl_0,\shl_1)
\to \rhom[{\WW[\stx]}](\shl_0,\shl_0)\simeq \cor_{\Lambda_0}
\eneqn
defines the morphism
\eq\label{eq:dualmorph0}
&&
\rhom[{\WW[\stx]}](\shl_0,\shl_1)\to
\rhom[\cor_{\Lambda_0}](\rhom[{\WW[\stx]}](\shl_1,\shl_0),\cor_{\Lambda_0}).
\eneq

To conclude, remark that for an object $F\in \RD^\Rb(\cor_\stx)$
supported by $\Lambda_0$, we have
\eqn
\RD_\stx F\simeq \rhom[\cor_{\Lambda_0}](F,\cor_{\Lambda_0})\,[\dim_\C\stx].
\eneqn

\vspace{0.2cm}
\noindent
(ii) Let us prove that \eqref{eq:dualmorph0} is an isomorphism
by adapting the proof of \cite[Prop.~5.1]{K1}.

Since this is a local problem, we may assume that 
$\stx=T^*X$, $X$ is an open subset of $\C^d$, $\Lambda_0=T^*_XX$ and $\shl_0=\sho_X^\tau$. 
Since $\rhom[{\WW[\stx]}](\sho_X^\tau,\shl_1)$ is 
constructible, we are reduced 
to prove the
isomorphisms for each $x\in X$
\eq\label{eq:dualmorph2}
&&H^j(\rhom[{\WW[\stx]}](\shl_1,\sho_X^\tau))_x\simeq 
\indlim[U\ni x](H^j\rsect_{c}(U;\rhom[{\WW[\stx]}](\sho_X^\tau,\shl_1))\,[2d])^\star
\eneq
where $U$ ranges over the family of Stein open neighborhoods of $x$ and 
${}^\star$ denotes the duality functor in the category of $\cor$-vector spaces. 

We choose a finite free resolution of $\shl_1$ on a neighborhood of 
$x$:
\eqn
&&0\to \shw_{T^*X}^{N_r}\to[\cdot P_{r-1}]\cdots\to[\cdot P_{0}] \shw_{T^*X}^{N_0}\to \shl_1\to 0
\text{ for some $r\geq 0$}.
\eneqn
For a sufficiently small 
holomorphically convex compact neighborhood  $K$ of $x$, the object
$\rsect(K;\rhom[\shw_{T^*X}](\shl_1,\sho_X^\tau))$ is represented by 
the complex
\eqn\label{eq:Ecomplex}
E^\bullet(K)\eqdot&& 0\to (\sect(K;\sho_X^\tau))^{N_0}\to[P_0\cdot]\cdots
\to[P_{r-1}\cdot](\sect(K;\sho_X^\tau))^{N_r}\to 0,
\eneqn
where $(\sect(K;\sho_{X}^\tau))^{N_0}$ stands in degree $0$.
Since
\eqn
&&\rhom[\shw_{T^*X}](\sho_X^\tau,\shw_{T^*X}\,[d])\simeq\Omega_X^\tau,
\eneqn
the object $\rhom[\shw](\sho_X^\tau,\shl_1\,[d])$
is represented by the complex
\eqn
&&0\to(\Omega_X^\tau)^{N_r}\to[\cdot P_{r-1}]\cdots
     \to[\cdot P_{0}](\Omega_X^\tau)^{N_0}\to 0,
\eneqn
where $(\Omega_X^\tau)^{N_0}$ stands in degree $0$. Hence, for a sufficiently small Stein open
neighborhood  $U$ of $x$,
$\rsect_c(U;\rhom[\shw](\sho_X^\tau,\shl_1\,[2d]))$
is represented by the complex 
\eqn\label{eq:Fcomplex}
F^\bullet_c(U)\eqdot&& 0\to (H^d_c(U;\Omega_X^\tau))^{N_r}
\to[\cdot P_{r-1}]\cdots
\to[\cdot P_{0}](H^d_c(U;\Omega_X^\tau))^{N_0}\to 0
\eneqn
where $(H^d_c(U;\Omega_{X}^\tau))^{N_0}$ stands in degree $0$.

Let $U_n$ be the open ball of $X$ centered at $x$ with radius $1/n$. 
By the hypothesis, all morphisms 
\eqn
&& F_c^\bullet(U_{p})\to F_c^\bullet(U_{n})
\eneqn
are quasi-isomorphisms for $p\geq n\gg 0$, 
and the cohomologies are finite-dimensional over $\cor$. 
Therefore, the hypotheses of Theorem~\ref{th:functan3bis} are satisfied by 
Lemma~\ref{le:vanishshotau1} and we get
\eqn
&&\inddlim[n]\Hom[\cor](F_c^\bullet(U_{n}),\cor)\simeq
\inddlim[n](F_c^\bullet(U_{n}))^\vee \text{ in }\II[{\RK^\Rb(\mdfn[\cor])}].
\eneqn
Applying Proposition~\ref{pro:serreOtau}, we
obtain 
\eqn
H^j(\rhom[\shw](\shl_1,\sho_X^\tau))_x
&\simeq&\indlim[n]H^j(E^\bullet(\overline U_{n}))\\
&\simeq&\indlim[n]H^j\Hom[\cor](F_c^\bullet(U_{n}),\cor)\\
&\simeq&\indlim[n](H^j\rsect_{c}(U;\rhom[{\WW[\stx]}](\sho_X^\tau,\shl_1))\,[2d])^\star.
\eneqn
\end{proof}

\begin{corollary}\label{cor:duality1}
In the situation of {\rm Theorem~\ref{th:duality1}~(ii)}, assume moreover
that $\stx$ is compact of complex dimension $2n$. 
Then the $\cor$-vector spaces 
$\Ext[{\WW[\stx]}]j(\shl_1,\shl_0)$ and  $\Ext[{\WW[\stx]}]{2n-j}(\shl_0,\shl_1)$ 
are finite-dimensional and dual to each other.
\end{corollary}


\section{Statement of the main theorem}

For two subsets $V$ and $S$ of the real manifold $X$, the normal cone 
$C(S,V)$ is well defined in $TX$. If $V$ is a smooth and closed
submanifold of $X$, one denotes by $C_V(S)$ the image of $C(S,V)$  
in the normal bundle $T_VX$.

Consider now a complex symplectic manifold $(\stx,\symplecto)$. 
The $2$-form $\symplecto$ gives the Hamiltonian 
 isomorphism $H$ from the cotangent
to the tangent bundle to $\stx$:
\eq\label{eq:sympiso1}
&& H\cl T^*\stx\isoto T\stx, \quad 
\langle\theta,v\rangle=\symplecto(v,H(\theta)), \quad
v\in T\stx,\,\theta\in T^*\stx.
\eneq
For a smooth Lagrangian
submanifold $\Lambda$ of $\stx$, 
the isomorphism \eqref{eq:sympiso1} induces an
isomorphism between the normal bundle to $\Lambda$ in $\stx$ and its
cotangent bundle:
\eq\label{eq:sympiso2}
&&T_\Lambda\stx\simeq T^*\Lambda.
\eneq
Let us recall a few notations and conventions
(see \cite{K-S1}).
For a complex manifold $X$ and a complex analytic
subvariety  $Z$ of $X$, one denotes by $Z_\reg$ the smooth locus of
$Z$, a complex submanifold  of $X$.
For a holomorphic $p$-form $\theta$ on $X$, one says that $\theta$
vanishes on $Z$ and one writes $\theta\vert_Z=0$ if $\theta\vert_{Z_\reg}=0$.

\begin{proposition}\label{pro:normalcone}
Let $\stx$ be a complex symplectic manifold and  
let $\Lambda_0$ and $\Lambda_1$ be two closed complex analytic
isotropic subvarieties of $\stx$. Then,
after identifying $T\stx$ and  $T^*\stx$ by \eqref{eq:sympiso1}, the
normal cone $C(\Lambda_0,\Lambda_1)$ is  a complex analytic $\C^\times$-conic 
isotropic subvariety of $T^*\stx$. 
\end{proposition}
Note that the same result holds for real analytic symplectic manifolds,
replacing ``complex analytic variety'' with ``subanalytic subset'' and 
``$\C^\times$-conic'' with ``$\R^+$-conic''.

First we need two lemmas.

\begin{lemma}
Let $X$ be a complex manifold and $\theta$ a $p$-form on $X$.
Let $Z\subset Y$ be closed subvarieties of $X$.
If $\theta\vert_{Y}=0$, then $\theta\vert_Z=0$.
\end{lemma}
\begin{proof}
By Whitney's theorem, we can find an open dense subset
$Z'$ of $Z_\reg$ such that 
\eqn
&&\left\{\parbox{330pt}{
for any
sequence $\{y_n\}_n$ in $Y_\reg$ such that it converges to a point $z\in Z'$
and $\{T_{y_n}Y\}_n$ converges to a linear subspace $\tau\subset T_zX$,
$\tau$ contains $T_zZ'$.}\quad\right.
\eneqn
Since $\theta$ vanishes on $T_{y_n}Y$, it vanishes also on 
$\tau$ and hence on $T_zZ'$.
\end{proof}

\begin{lemma}\label{le:normalcone}
Let $X$ be a complex manifold, $Y$ a closed complex subvariety of $X$
and $f\cl X\to \C$ a holomorphic function. Set $Z\eqdot\opb{f}(0)$,
$Y'\eqdot \overline{(Y\setminus Z)}\cap Z$. Consider a $p$-form
$\eta$, a $(p-1)$-form $\theta$ on $X$ and set
\eqn
&&\omega=df\wedge \theta +f\eta.
\eneqn
Assume that  $\omega\vert_Y=0$. Then $\theta\vert_{Y'}=0$ and 
 $\eta\vert_{Y'}=0$.
\end{lemma} 

\begin{proof}  
We may assume that $Y=\overline{(Y\setminus Z)}$.

Using Hironaka's desingularization theorem, we may find a smooth manifold
$\tw Y$ and a proper morphism $p\cl\tw Y\to Y$ such that,
in a neighborhood of each point of $\tw Y$, 
$p^*f$ may be written in a local coordinate system $(y_1,\dots,y_n)$ as a
product $\prod_{i=1}^ny_i^{a_i}$, where the $a_i$'s are non-negative integers. 
Let $\tw Z=\opb{p}(Z)$. Then 
$\overline{(\tw Y\setminus \tw Z)}\cap\tw Z\to Y'$ is proper and surjective.

Hence, we may assume from the beginning that $Y$ is smooth and 
then $Y=X$.
Moreover, 
it  is enough to prove the result at generic points of $Y'$. Hence, we 
may assume, setting  $(y_1,\dots,y_n)=(t,x)$ ($x=(y_2,\dots,y_n)$), 
that $f(t,x)=t^a$ for some $a>0$. 
Write
\eqn
&&\theta=t\theta_0+dt\wedge\theta_1+\theta_2,
\eneqn
where $\theta_1$ and $\theta_2$ depend only on $x$ and $dx$.

By the hypothesis, 
\eqn
&&0=df\wedge \theta+ f\eta
=at^adt\wedge\theta_0+at^{a-1}dt\wedge\theta_2+ t^a\eta.
\eneqn
It follows that $\theta_2$ is identically zero.
Hence:
\eqn
&&\theta=t\theta_0+dt\wedge\theta_1,\quad\eta=-a dt\wedge\theta_0.
\eneqn
Therefore, $ \theta\vert_{t=0}=\eta\vert_{t=0}=0$.
\end{proof}

\renewcommand{\proofname}{Proof of Proposition~\ref{pro:normalcone}}
\begin{proof}
Recall that $\stx^a$ denotes the complex manifold $\stx$ endowed with
the symplectic form $-\symplecto$ and that 
$\Delta$ denotes the diagonal of $\stx\times\stx^a$.
Using the isomorphisms
\eqn
&&T\stx\simeq T_\Delta(\stx\times\stx^a),\\
&&C(\Lambda_0,\Lambda_1)\simeq C(\Delta,\Lambda_0\times\Lambda_1^a),
\eneqn
(the first isomorphism is associated with the first projection 
on $\stx\times\stx^a$)
we are reduced to prove the result when $\Lambda_0$ is smooth.

Let $(x,u)$ be a local symplectic coordinate system on $X$ such that
\eqn
 \Lambda_0=\set{(x;u)}{u=0},
&&\symplecto=\sum_{i=1}^n du_i\wedge dx_i.
\eneqn
Consider the normal deformation  $\tw\stx_{\Lambda_0}$ of 
$\stx$ along $\Lambda_0$. 
Recall that we have a diagram
\eqn
\xymatrix{
{\stx} &{\tw\stx_{\Lambda_0}}\ar[l]_-p\ar[r]^-t&{\C}\\
{\Lambda_1}\ar@{^{(}->}[u]&{\opb{p}\Lambda_1}\ar@{^{(}->}[u]\ar[l]_-p&
}\eneqn
such that, denoting by $(x;\xi,t)$ the coordinates on
$\tw\stx_{\Lambda_0}$,
\eqn
&&p(x;\xi,t)=(x;t\xi),\\
&&T_{\Lambda_0}\stx\simeq \set{(x;\xi,t)\in \tw\stx_{\Lambda_0}}{t=0},\\
&&C_{\Lambda_0}(\Lambda_1)\simeq\overline{\opb{p}(\Lambda_1)\setminus\opb{t}(0)}\cap\opb{t}(0).
\eneqn
Clearly, $C_{\Lambda_0}(\Lambda_1)$ is a complex analytic
variety.
Moreover,  
\eqn
&&p^*\symplecto=\sum_{i=1}^n d(t\xi_i)\wedge dx_i
=dt\wedge\bigl(\sum_{i=1}^n \xi_i\wedge dx_i\bigr)
+t\sum_{i=1}^nd\xi_i\wedge dx_i.
\eneqn
Since $p^*\symplecto$ vanishes on $(\opb{p}\Lambda_1)_\reg$, $\sum_{i=1}^nd\xi_i\wedge dx_i$
vanishes on $C_{\Lambda_0}(\Lambda_1)$ by Lemma~\ref{le:normalcone}.
\end{proof}
\renewcommand{\proofname}{Proof}

\begin{theorem}\label{th:main1}
Let $\stx$ be a complex symplectic manifold and let 
$\shl_i$ \lp$i=0,1$\rp\, be two objects of 
$\RD^\Rb_{\rh}(\WW[\stx])$ supported by smooth Lagrangian 
submanifolds  $\Lambda_i$. Then
\bnum
\item
the object $\rhom[{\WW[\stx]}](\shl_1,\shl_0)$ 
belongs to  $\RD^\Rb_\Cc(\cor_\stx)$ and its microsupport is contained 
in the normal cone $C(\Lambda_0,\Lambda_1)$,
\item
the natural morphism
\eqn
&&\rhom[{\WW[\stx]}](\shl_1,\shl_0)\to
\RD_\stx\bigl(\rhom[{\WW[\stx]}](\shl_0,\shl_1\,[\dim_\C\stx])\bigr)
\eneqn
is an isomorphism.
\enum
\end{theorem}
The proof of (i)  will be given in \S~\ref{section:proof1}
and (ii) is a particular case of Theorem~\ref{th:duality1}.

\begin{conjecture}\label{conj:mainconj}
Theorem~\ref{th:main1} remains true without assuming that 
the $\Lambda_i$'s are smooth.
\end{conjecture}

Remark  that the analogous of Conjecture~\ref{conj:mainconj} for
complex contact manifolds is true over the field $\C$, as we shall see
in \S~\ref{section:contact}.

\begin{corollary}\label{cor:main1}
Let  $\shl_0$ and $\shl_1$ be two regular holonomic 
${\WW[\stx]}$-modules supported by smooth Lagrangian submanifolds.
Then the object 
$\rhom[{\WW[\stx]}](\shl_1,\shl_0)$ of $\RD^\Rb_\Cc(\cor_\stx)$
is perverse.
\end{corollary}

\begin{proof}
Since the
problem is local,  we may assume that $\stx=T^*X$, $\Lambda_0=T^*_XX$ and
$\shl_0=\sho_X^\tau$. 

By Theorem~\ref{th:main1}~(ii), it is enough to check that if $Y$ is a
locally closed smooth submanifold of $X$ of codimension $l$, then  
$H^j(\rsect_Y(\rhom[\shw_{T^*X}](\shl_1,\shl_0))\vert_Y$ vanishes for $j<l$. This
follows from Lemma~\ref{le:vanisectYotau}.
\end{proof}

\begin{remark}\label{rem:BF}
It would be
interesting to compare $\rhom[{\WW[\stx]}](\shl_1,\shl_0)$
with the complexes obtained in
\cite{B-F}. 
\end{remark}

\section{Proof of Theorem~\ref{th:main1}}\label{section:proof1}

In this section $X$ denotes a complex manifold. As usual, $\sho_X$ is
the structure  sheaf and $\shd_X$ is the sheaf of rings of 
(finite--order) differential operators. For a coherent 
$\shd_X$-module $\shm$, we denote by $\chv(\shm)$ its characteristic
variety. This notion extends to the case where 
 $\shm$ is a countable union of coherent sub-$\shd_X$-modules. In this 
 case, one sets
\eqn
&&\chv(\shm)=\overline{\bigcup_{\shn\subset\shm}\chv(\shn)}
\eneqn
where $\shn$ ranges over the family of coherent $\shd_X$-submodules of 
$\shm$, and, for a subset $S$ of $T^*X$,  $\overline S$ means the closure 
of $S$. 

\begin{lemma}\label{le:SScar1}
Let $\shm_0$ be a coherent  $\shd_X$-module. Then 
\eq\label{eq:prop1}
&&SS(\rhom[\shd_X](\shm_0,\sho_X^\tau(0)))\subset\chv(\shm_0).
\eneq
\end{lemma}
\begin{proof}
Apply \cite[Ch.3~Th.~3.2.1]{S}.
\end{proof}

\begin{lemma}\label{le:SScar2}
Let $\shm$ be a coherent $\Dotau$-module.  Then 
\eq\label{eq:prop3}
&&SS(\rhom[{\Dotau}](\shm,\sho_X^\tau(0)))\subset\chv(\shm).
\eneq
\end{lemma}
\begin{proof}
Let $\shm_0\subset\shm$ be a coherent $\shd_X$-submodule which generates
$\shm$. Then 
\eqn
&&\chv(\Dotau\tens[\shd_X]\shm_0)=\chv(\shm_0)\subset\chv(\shm).
\eneqn

Consider the exact sequence of coherent $\Dotau$-modules
\eq\label{eq:esqshm0}
&& 0\to\shn\to \Dotau\tens[\shd_X]\shm_0\to\shm\to 0.
\eneq
Applying the functor 
$\rhom[{\Dotau}](\scbul,\sho_X^\tau)$ 
to the exact sequence \eqref{eq:esqshm0}, we get 
a distinguished triangle $G'\to G\to G''\to[+1]$\,. 
Note that $G\simeq \rhom[\shd_X](\shm_0,\sho_X^\tau)$, since 
$\Dotau$ is flat over $\shd_X$. 

Let $\theta=(x_0;p_0)\in T^*X$ with $\theta\notin \chv(\shm)$ 
and let $\psi$ be a
real function on $X$ such that $\psi(x_0)=0$ and $d\psi(x_0)=p_0$. 
Consider the distinguished triangle
\eq\label{eq:dt1}
&&(\rsect_{\psi\geq 0}(G'))_{x_0}\to (\rsect_{\psi\geq 0}(G))_{x_0}\to 
(\rsect_{\psi\geq 0}(G''))_{x_0}\to[+1]\, .
\eneq
By Lemma \ref{le:SScar1}, we have:
\eqn
&& H^j((\rsect_{\psi\geq 0}(G))_{x_0})\simeq 0\mbox{ for all }j\in\Z.
\eneqn
The objects of the distinguished triangle \eqref{eq:dt1} 
are concentrated in degree
$\geq 0$. Therefore, 
 $H^j((\rsect_{\psi\geq 0}(G'))_{x_0})\simeq 0$ for $j\leq 0$.

Since
$\chv(\shn)\subset\chv(\Dotau\tens[\shd_X]\shm_0)$,
we get
$H^j((\rsect_{\psi\geq 0}(G''))_{x_0})\simeq 0$ for $j\leq 0$.
By repeating this argument, we deduce that 
$H^j((\rsect_{\psi\geq 0}(G'))_{x_0})\simeq H^{j-1}((\rsect_{\psi\geq 0}(G''))_{x_0})\simeq 0$
for all $j\in\Z$.
\end{proof}

Note that the statement of Theorem~\ref{th:main1} is local and 
invariant by quantized symplectic transformation.
{}From now on, we denote by $(x;u)$  a local symplectic coordinate system 
on $\stx$ such that 
\eqn
&&X\eqdot\Lambda_0=\set{(x;u)\in \stx}{u=0}.
\eneqn
We denote by $(x;\xi)$ the associated homogeneous symplectic coordinates on
$T^*X$. The differential operator $\partial_{x_i}$ on $X$ 
has order $1$ and principal
symbol $\xi_i$. The monomorphism \eqref{eq:DinW} extends as a
monomorphism of rings
\eqn
&&\Dttau\hookrightarrow \W[T^*X].
\eneqn
Note that the total symbol of the operator 
$\partial_{x_i}$ of $\W[T^*X]$ is $u_i\cdot\tau$. 

We may assume that
there exists a holomorphic function $\phi\cl X\to\C$ such that 
\eq\label{eq:lambda1}
&&\Lambda_1=\set{(x;u)\in \stx}{u=\grad\,\phi(x)}.
\eneq
Here, 
$\grad\,\phi=(\phi'_1,\dots,\phi'_n)$  and  
$\phi'_i=\frac{\partial\phi}{\partial{x_i}}$.

If $\Lambda_0=\Lambda_1$, Theorem \ref{th:main1} is
immediate. 
We shall assume that  $\Lambda_0\neq\Lambda_1$ and thus that $\phi$ is not
a constant function.
We may assume that
\eq\label{eq:mod}
&&\left\{\parbox{300 pt}{
$\shl_0=\OO^\tau_X$,\\[3pt]
$\shl_1=\W[T^*X]/\shi_1$, 
$\shi_1$ being the ideal generated by 
$\{\partial_{x_i}-\phi'_i\tau\}_{i=1,\dots,n}$.
}\right. 
\eneq
To $\phi\cl X\to\C$ are associated the maps
\eqn
&&T^*X\from[\phi_d]X\times_\C T^*\C\to[\phi_\pi]T^*\C,
\eneqn
and the $(\shd_X,\shd_\C)$-bimodule $\shd_{X\to[\phi]\C}$.
Let $t$ be the coordinate on $\C$. By
identifying $\partial_t$ and $\tau$, we regard
$\shd_{X\to[\phi]\C}$ as a $\shd_X[\tau]$-module.
We set
\eqn
&&V=\overline{\im \phi_d}, 
\mbox{ the closure of $\phi_d(T^*\C\times_\C X)$}.
\eneqn

\begin{lemma}\label{le:prop3}
Regarding $\shd_{X\to[\phi]\C}$ as a $\shd_X$-module, 
one has $\chv(\shd_{X\to[\phi]\C})\subset V$.
\end{lemma}
Note that $\shd_{X\to[\phi]\C}$ is not a coherent $\shd_X$-module in
general.

\begin{proof}
Let $\Gamma_\phi$ be the graph of $\phi$ in $X\times\C$ 
and let $\shb_\phi$ be the coherent left $\shd_{X\times\C}$-module associated
 to this submanifold. 
(Note that  $\chv(\shb_\phi)=\Lambda_\phi$, the 
conormal bundle to $\Gamma_\phi$ in $T^*(X\times\C)$.)
We shall identify $\Gamma_\phi$ with $X$ by the first projection  on $X\times\C$
and the left $\shd_X$ modules $\shb_\phi$ with $\shd_{X\to[\phi]\C}$.

Denote by $\delta(t-\phi)$ the canonical generator of $\shb_\phi$. Then
\eqn
&&\shd_{X\to[\phi]\C}= \shd_{X\times\C}\cdot\delta(t-\phi)=
\sum_{n\in\N}\shd_X\partial_t^n\cdot\delta(t-\phi).
\eneqn
Define
\eqn
&&\shn_\phi\eqdot
\sum_{n\in\N}\shd_X(t\partial_t)^n\cdot\delta(t-\phi).
\eneqn
By \cite[Th.~6.8]{K3}, 
the $\shd_X$-module $\shn_\phi$ is coherent and its characteristic
variety is contained in $V$.  All $\shd_X$-submodules $\partial_t^n\shn_\phi$ of
$\shd_{X\to[\phi]\C}$ are isomorphic since
$\partial_t$ is injective on $\shb_\phi$. Then the result follows from 
\eqn
&&\shd_{X\to[\phi]\C}=\sum_{n\in\N}\partial_t^n\shn_\phi.
\eneqn
\end{proof}

We consider the left ideal and the modules:
\eq\label{eq:shj1}
&&\left\{  \parbox{300 pt}{
$\shi \eqdot$ the left ideal of $\Dttau$ 
generated  by $\{\partial_{x_i}-\phi'_i\tau\}_{i=1,\dots,n}$,\\
$\shm =\Dotau/(\shi\cap\Dotau)$,\\
$\shm_0= \shd_X/(\shi\cap\shd_X)$.
}\right. 
\eneq

\begin{lemma}\label{le:SScar3}
\bnum
\item
$\shm$ is $\Dotau$-coherent,
\item
$\shm_0$ is $\shd_X$-coherent,
\item
we have an isomorphism of $\shd_X$-modules 
$\shd_X\cdot\delta(t-\phi)\simeq\shm_0$, where 
$\shd_X\cdot\delta(t-\phi)$ is the $\shd_X$-submodule 
of $\shd_{X\to[\phi]\C}$ generated by $\delta(t-\phi)$,
\item
$\chv(\shm)= \chv(\shm_0)\subset V$.
\enum
\end{lemma}
\begin{proof}
(i) follows from Lemma \ref{le:dtau2} (ii).

\noindent
(ii) follows from Lemma \ref{le:dtau1} (iii).

\noindent
(iii) Clearly, $\shd_X\cdot\delta(t-\phi)\simeq
\shd_X/\shd_X\cap\shi$. 

\noindent
(iv) (a) By (iii) and Lemma~\ref{le:prop3} we get the inclusion
$\chv(\shm_0)\subset V$. 

\noindent
(iv) (b) Since $\shm_0\subset\shm$, the inclusion 
$\chv(\shm_0)\subset\chv(\shm)$ is obvious.

\noindent
(iv) (c) Denote by $u$ the canonical generator of $\shm$. Then 
$\shm=\bigcup_{n\leq 0}\shd_X\tau^nu$. Since there are 
epimorphisms $\shm_0\epito\shd_X\tau^nu$, the inclusion 
$\chv(\shm)\subset\chv(\shm_0)$ follows.
\end{proof}

Set
\eq\label{eq:F0}
&&F_0=\rhom[{\Dotau}](\shm,\sho_X^\tau(0)),
\eneq
where the module $\shm$ is defined in \eqref{eq:shj1}.

\begin{lemma}\label{le:mainthL0a}
$SS(F_0)$ is a closed  $\C^\times$-conic complex analytic Lagrangian subset of
$T^*X$ contained in $C(\Lambda_0,\Lambda_1)$.
\end{lemma}
\begin{proof}
By Lemmas \ref{le:SScar2} and \ref{le:SScar3},
$SS(F_0)\subset V\cap (\opb{\pi}(\Lambda_0\cap \Lambda_1))$, and one
immediately checks that
\eq\label{eq:Vcap=C}
&&  V\cap \opb{\pi}(\Lambda_0\cap \Lambda_1)=C(\Lambda_0,\Lambda_1).
\eneq

Since $SS(F_0)$ is involutive by \cite{K-S1} and is contained in a 
$\C^\times$-conic analytic isotropic subset by Proposition~\ref{pro:normalcone},
it is a closed  $\C^\times$-conic 
complex analytic Lagrangian subset of $T^*\Lambda_0$ by 
\cite[Prop.~8.3.13]{K-S1}.
\end{proof}

\begin{lemma}\label{le:mainthL0b}
Let $F_0$ be as in \eqref{eq:F0}. Then for each $x\in X$ and each
$j\in\Z$, the $\coro$-module $H^j(F_0)_x$ is finitely generated.
\end{lemma}

\begin{proof}
Let $x_0\in X$ and choose a local coordinate system around $x_0$.
Denote by $B(x_0;\epsilon)$ the closed ball of center $x_0$ and radius
$\epsilon$.  By a result of \cite{K1} (see also \cite[Prop.~8.3.12]{K-S1}), 
we deduce from Lemma~\ref{le:mainthL0a} that the natural morphisms
\eqn
&&\rsect(B(x_0;\epsilon_1);F_0)\to \rsect(B(x_0;\epsilon_0);F_0)
\eneqn
are isomorphisms for  $0< \epsilon_0\leq \epsilon_1\ll 1$.

We represent $F_0$ by  a complex:
\eq\label{eq:cpxF0}
&& 0\to (\sho^\tau(0))^{N_0}\to[d_0]\cdots\to (\sho^\tau(0))^{N_n}\to 0,
\eneq
where the differentials are $\coro$-linear.
It follows from  Lemma~\ref{le:functan2} and Theorem~\ref{th:functan1}
that the cohomology objects $H^j(F_0)_{x_0}$ are finitely generated.
\end{proof}


\renewcommand{\proofname}{End of the proof of Theorem~\ref{th:main1}}
\begin{proof}  
As already mentioned, part (ii) is a particular case of Theorem~\ref{th:duality1}.

Let us prove part (i).
By ``d{\'e}vissage'' we may assume that $\shl_0$ and $\shl_1$ are
concentrated in degree $0$. We may also assume that 
$\shl_0$ and $\shl_1$ are as in \eqref{eq:mod}.

Let $F_0$ be as in \eqref{eq:F0}, and
set 
$F=\rhom[\shw_{T^*X}](\shl_1,\sho_X^\tau)$.
Since $\shl_1\simeq\shw_{T^*X}\tens_{\Dotau}\shm$, we obtain
$$F\simeq F_0\tens[\coro]\,\cor$$
by Lemmas~\ref{le:dtau2} and~\ref{le:dtau3}.

Hence we have $SS(F)\subset SS(F_0)$,
and the weak constructibility  as well as 
the estimate of $SS(F)$ follows by 
Lemma \ref{le:SScar2} and Lemma~\ref{le:mainthL0a}. 
The finiteness result follows from Lemma~\ref{le:mainthL0b}.
\end{proof}
\renewcommand{\proofname}{Proof}

\section{Serre functors on contact and symplectic manifolds}
\label{section:contact}

In the definition below, $\kor$ is a field and ${}^\star$ denotes the
duality functor for $\kor$-vector spaces. 

\begin{definition}\label{def:BKcat}
Consider a $\kor$-triangulated category $\sht$.
\bnum
\item
The category $\sht$ is $\Ext[]{}$-finite if for any $L_0,L_1\in\sht$, 
$\Ext[\sht]j(L_1,L_0)$ is finite-dimensional over $\kor$ for all
$j\in\Z$ and is zero for $\vert j\vert\gg 0$.
\item 
Assume that $\sht$ is $\Ext[]{}$-finite. A Serre functor 
(see~\cite{B-K})
$S$ on $\sht$ is  an equivalence of $\kor$-triangulated categories 
$S\cl\sht\to\sht$ satisfying 
\eqn
&&\bigl(\Hom[\sht](L_1,L_0)\bigr)^\star\simeq\Hom[\sht](L_0,S(L_1))
\eneqn
functorially in $L_0,L_1\in\sht$.  
\item 
If moreover there exists an integer $d$ such that $S$ is isomorphic to
the shift by $d$, then one says that $\sht$ is a $\kor$-triangulated
Calabi-Yau category of dimension $d$.
\enum
\end{definition}

Let $\sty$ be a complex contact manifold. 
The algebroid stack $\EE[\sty]$
of microdifferential operators on $\sty$ has been constructed in \cite{K2} and
the triangulated categories
$\RD^\Rb_{\coh}(\EE[\sty])$, $\RD^\Rb_{\hol}(\EE[\sty])$ and
$\RD^\Rb_{\rh}(\EE[\sty])$ are naturally defined.

\begin{theorem}\label{th:main3}
For a complex contact manifold  $\sty$, we have 
\bnum
\item
for $\shm$ and $\shn$ in 
$\RD^\Rb_{\rh}(\EE[\sty])$, the object 
$F=\rhom[{\EE[\sty]}](\shm,\shn)$ belongs to $\RD^\Rb_{\Cc}(\C_\sty)$,
\item
if $\sty$ is compact, then
$\RD^\Rb_{\rh}(\EE[\sty])\,$ is a Calabi-Yau
$\C$-triangulated category of dimension $\dim_\C\sty-1$. 
\enum
\end{theorem}
\renewcommand{\proofname}{\em Sketch of proof}
\begin{proof}
(i) is well-known and follows from \cite{K-K} (see \cite{W} for
further developments). The idea of the proof is as follows.
The assertion being local and
invariant by  quantized contact
transformations, we may assume that
$\sty$ is an open subset of the projective cotangent bundle $P^*Y$ to a complex
manifold. Then, using the diagonal procedure, we reduce to the case 
$F=\rhom[\shd_Y](\shm,\shc_{Z\vert Y})$, where $\shm$ is a regular
holonomic $\shd_Y$-module and $\shc_{Z\vert Y}$ is the
$\she_{P^*Y}$-module 
associated to a complex hypersurface $Z$ of $Y$. 

\noindent
(ii) follows from (i) as in the proof of Theorem~\ref{th:duality1}.
\end{proof} 
\renewcommand{\proofname}{Proof} 

\begin{remark}
(i) If Conjecture~\ref{conj:mainconj} is true, that is, 
if Theorem~\ref{th:main1}
holds for any Lagrangian varieties, then, for any compact complex
symplectic manifold $\stx$, 
the $\cor$-triangulated category 
$\RD^\Rb_{\rh}(\WW[\stx])$ is a Calabi-Yau
$\cor$-triangulated category of dimension $\dim_\C\stx$. Note that 
this result is true when replacing the notion of regular holonomic
module by the notion of good 
modules, i.e., 
coherent modules admitting globally defined 
$\WW[\stx](0)$-submodules which generate them. This follows from
a theorem of Schapira-Schneiders to appear.

\noindent
(ii) Note that Proposition~1.4.8 of \cite{K-K} is not true, but this
proposition is not used in loc.\ cit. Indeed,
for  a compact complex manifold $X$, if $\sht\eqdot
\RD^\Rb_{\rh}(\shd_X)$ denotes the 
full triangulated subcategory of  $\RD^\Rb(\shd_X)$ consisting 
of objects with regular holonomic cohomologies, it is well known that 
the duality functor is not a Serre functor on $\sht$. 

\noindent
(iii) It may be interesting to notice that the duality functor 
is not a Serre functor in the setting of regular holonomic $\shd$-modules, 
but is a Serre functor in
the microlocal setting, that is, for regular holonomic $\she$-modules. This may be
compared to  \cite[Prop.~8.4.14]{K-S1}.
\end{remark}

\vspace*{1cm}
\noindent
\parbox[t]{17em}
{\scriptsize{
\noindent
Masaki Kashiwara\\
Research Institute for Mathematical Sciences \\      
Kyoto University \\                       
Kyoto, 606--8502, Japan\\         
e-mail: masaki@kurims.kyoto-u.ac.jp
}}
\quad
\parbox[t]{14em}
{\scriptsize{
Pierre Schapira\\
Institut de Math{\'e}matiques\\
Universit{\'e} Pierre et Marie Curie\\
175, rue du Chevaleret,
75013 Paris, France\\
e-mail: schapira@math.jussieu.fr\\
http://www.math.jussieu.fr/$\sim$schapira/}}


\begin{thebibliography}{15}


\bibitem{B-F} K.~Behrend and B.~Fantechi,
{\em  Gerstenhaber and Batalin-Vilkovisky structures on 
Lagrangian intersections,}
to appear in ``Arithmetic and Geometry, Manin Festschrift''.\\
http://www.math.nyu.edu/{\~{}}tschinke/.manin/manin-index.html

\bibitem{B-K} A.~I.~Bondal and M.~Kapranov, 
{\em Representable functors, Serre functors and mutations,} 
(English translation) Math USSR Izv., {\bf 35} p.~519--541 (1990).

\bibitem{D-P} A.~D'Agnolo and P.~Polesello,
{\em Deformation-quantization of complex involutive submanifolds,}
in:  Noncommutative Geometry and Physics,
World Scientific  p.~127-137 (2005).

\bibitem{D-S} A.~D'Agnolo and P.~Schapira,
{\em Quantization of complex Lagrangian submanifolds,}
ArXiv:math. AG/0506064.

\bibitem{Gr} A.~Grothendieck,
{\em Produits Tensoriels Topologiques et Espaces Nucl{\'e}aires,} 
Mem. Am. Math. Soc. {\bf 16} (1955); erratum,
Ann. Inst. Fourier {\bf 6} p.~117--120 (1955/56).

\bibitem{Ho} C.~Houzel,
{\em Espaces analytiques relatifs et th{\'e}or{\`e}mes de finitude,}
Math. Annalen {\bf 205} p.~13--54 (1973).

\bibitem{K1} M.~Kashiwara,
{\em On the maximally overdetermined systems of linear differential equations,}
Publ. RIMS, Kyoto Univ. {\bf 10} p.~563--579 (1975).

\bibitem{K2} \bysame,
{\em Quantization of contact manifolds,} 
Publ. RIMS, Kyoto Univ. {\bf 32} p.~1-5 (1996).

\bibitem{K3} \bysame,
{\em D-modules and Microlocal Calculus,}
Translations of Mathematical Monographs, {\bf 217} American
Math. Soc. (2003).

\bibitem{K-K} M.~Kashiwara and T.~Kawai,
{\em On holonomic systems of microdifferential equations III,}
Publ. RIMS, Kyoto Univ. {\bf 17} p.~813--979 (1981).

\bibitem{K-O}
M.~Kashiwara and T.~Oshima,
{\em Systems of differential equations with regular singularities and
 their boundary value problems,}
Ann. of Math. {\bf 106} p.~145--200 (1977).


\bibitem{K-S1} M.~Kashiwara and P.~Schapira,
{\em Sheaves on Manifolds,}
Grundlehren der Math. Wiss. {\bf 292} Springer-Verlag (1990).

\bibitem{K-S2} \bysame,
{\em Categories and Sheaves,}
Grundlehren  der Math. Wiss. {\bf 332} Springer-Verlag (2005).

\bibitem{K-V} Kiel and J-L.~Verdier,
{\em Ein einfascher beweis der koh{\"a}renzsatzes von Grauert,}
Math. Annalen {\bf 195} p.~24--50 (1971).

\bibitem{Ko} M.~Kontsevich, 
{\em Deformation quantization of algebraic varieties}
Lett. Math. Phys. {\bf 56} p.~271--294, (2001).



\bibitem{P-S} P.~Polesello and P.~Schapira,
{\em Stacks of quantization-deformation modules over
  complex symplectic manifolds,}  Int. Math. Res. Notices {\bf 49} 
p.~2637--2664 (2004).


\bibitem{Pr-S} F.~Prosmans and J-P.~Schneiders,
   {\em   A topological reconstruction theorem for $D^{\infty}$-Modules},
     Duke Math. J. {\bf 102} p.~39-86 (2000).


\bibitem{S-K-K} M.~Sato, T.~Kawai, and M.~Kashiwara,
{\em Microfunctions and pseudo-differential equations,}
in Komatsu (ed.), {\em Hyperfunctions and pseudo-differential  equations,}
Proceedings Katata 1971, Lecture Notes in Math. Springer-Verlag
{\bf 287} p.~265--529 (1973).

\bibitem{S} P.~Schapira,
{\em Microdifferential Systems in the Complex Domain,}
Grundlehren der Math. Wiss. {\bf 269} Springer-Verlag (1985).

\bibitem{Scn0} J-P.~Schneiders,
{\em Quasi-abelian Categories and Sheaves,}  
Mem. Soc. Math. Fr. {\bf 76} (1999).

\bibitem{W} I.~Waschkies,
{\em Microlocal Riemann-Hilbert correspondence,}
Publ. RIMS, Kyoto Univ. {\bf 41} p.~35-72 (2005).

\end{thebibliography}
\end{document}